\documentclass{article}
\newtheorem{theorem}{Theorem}
\usepackage{graphicx}
\usepackage{color}
\usepackage{graphicx,epsfig}
\usepackage{epstopdf}
\usepackage{amsmath}
\usepackage{epstopdf}
\title{Finite element analysis of modified N-S equations coupled with energy transfer for Hybrid Nanofluid flow in complex domains}
\author{Sangita Dey, B.V. Rathish Kumar\\
Dept. of Mathematics and Statistics, IIT Kanpur, 
           ,Kanpur,
            208016, 
            U.P,India}
\date{March 2023}

\begin{document}

\maketitle

\begin{abstract}
A theoretical and computational finite element study of modified Navier-Stokes Equations coupled with energy conservation governing the flow and heat transfer in complex domains with hybrid nanofluid$(HNF)$ is carried out. The apriori error estimates providing the convergence analysis for the finite element scheme is derived in the $H^1$ norm.   A detailed parametric analysis on the physics related to flow and heat tranfer for different governing parameters such as: hybrid nano-particle volume fraction $(\phi),$ Rayleigh Number $(Ra)$,  Prandtl Number$(Pr)$,  and the length of heat source has been reported. The use of hybrid nano-particles leads to the enhancement of heat conduction both in L and H shaped domains. While a butterfly flow pattern with a multi-cellular structure is noticed in H-shaped domain, the flow in L-Shaped domain depicts a boundary layer feature with primary, secondary and tertiary co-rotating and counter rotating cells. 
\end{abstract}



\section{keyword}
Galerkin Method $\cdot$ Convection diffusion $\cdot$ Rayleigh number $ \cdot $ Prandtl number $\cdot$  a priori estimates $\cdot$  Hybrid nanofluid.
 
\section{Introduction}
\label{sec1}
Numerical solution through finite element method is a popular approach to study  partial differential equations (PDEs). The trustworthiness of the derived numerical approximation to the PDEs is established by the investigation of the numerical schemes' convergence. 
The derivation of the  apriori finite element error estimates provides a robust mathematical approach to establish such a convergence analysis. 
These apriori error analysis guarantees that the numerical solution is both bounded and converges with increased spatial resolutions to the exact solution. Applications for such a coupled nonlinear PDE system include solar collectors, ground water pollution transmission, storage of nuclear waste, crude oil production, and geothermal systems \cite{SRS91,CC89,Sm17}. Many such industrial applications are well known to suffer from the low thermal conductivity problem since several decades. Nanotechnology has attracted a lot of attention in an effort to overcome these constraints, and the study of nanofluid flow has become increasingly important recently. 
Preparation of nanofluid needs suspension of nanoparticles in the base fluids like water, oil, ethylene glycol etc. There are a few commonly used nanofluids like $Al_2O_3$/water, $Fe_2O_3$/water, Cu/water, Ag/water, $Ti_2O_3$/water etc. In \cite{MMA10} authors have concluded that the average Nusselt no. is maximum for $Al_2O_3$/water compared to Cu/water, Ag/water, $Ti_2O_3$/water nanofluid. Also several studies can be found on hybrid nanofluids where mixture of two or more different types of nanoparticles are used to get much better result. In this regards $Cu-Al_2O_3/$water is a commonly used hybrid nanofluid to serve various industrial problems. One such work can be found in \cite{YRMZ19} where effect of the magnetic fields on the hybrid nanofluid (HNF) Ag-Mgo/water is analyzed in a channel attached with a heat source and cooling sink. The effect of various input parameters, on heat transfer, like Hartmann number $(0\leq Ha \leq 60)$,
Reynolds number $(0\leq Re \leq 100)$, nano particles volume fractions $(0\leq \phi \leq 0.02)$ e.t.c are analyzed. L shaped geometry is one of the important geometry which has lots of practical significance. In their investigation of the impact of the position and size of the heat source on magneto-hydrodynamics, the authors $[4]$ take into account an L-shaped geometry. A thorough analysis of hydrodynamic convection in
$Al_2O_3/$water is conducted. 
In the paper\cite{MSR21} authors have considered an open ended L shaped domain filled up with $CuO-H_{2}O$ nano-fluid to analyze the influence of the partinent parameters like cavity aspect ratio $0.2 \leq AR \leq 0.6$, volume fraction of nanoparticles $(0 \leq \phi \leq 0.1)$, Rayleigh number $(10^3 \leq Ra \leq 10^6 )$ on heat transfer rate. In \cite{TMAI19} authors have chosen a new type of nanofluid $Ag-water$ nanofluid in an L shaped channel having an inclination $30^{0}$ with X-axis. The governing equations are solved by using finite volume method where an isothermal source was kept in the mid point of right inclined wall by keeping other walls adiabatic. The main conclusion of their experiment is that the effect of the aspect
ratio, nanoparticle volume percentage, Raynolds number, and Richardson number raises the mean Nusselt number and takes the entropy generation rate to a maximum. The authors in \cite{SVSC12} aim to investigate the effect of the location of a discrete heat source in an inclined open L shaped enclosure filled up with $Ag-water$ nanofluid. Here authors claim that the increment in the inclination angle increases the heat transfer. In \cite{BMA16}  the authors carry out a numerical experiment in an L shaped domain filled up with $Al_2O_3$/water where a uniform magnetic field is present. The effect of governing parameters like Rayleigh number, Hartmann number, nanoparticles' volume fraction, the cavity aspect ratio, on heat conduction, is presented.
\cite{TANM16} gives us the opportunity to see the effect of $Al_2O_3$/water nanofluid in a baffled L shaped enclosure on natural convection and heat transfer by finite volume method.\\ 
Yet another complex domain of industrial significance is the H-Shaped domain. Applications like as printing, HVAC, solar collectors, and other fields are
pertinent to the study of natural convection in these fields etc. \cite{ASA21}. This H-shaped enclosure is known for its effectiveness in heat treatment \cite{ASA21, FPM20,ARM18}. In \cite{ASA21} the authors have investigated the influence of fin and corrugated walls on double diffusive natural convection in a $(Cu, Al_2O_3)$/water filled H shaped cavity.\\

In this paper, we derive the apriori error estimates for the finite element approximation of the modified N-S equations coupled with energy conservation which are governing the fluid flow and  heat transfer in the presence of hybrid nano particles. Even though a plenty of literature can be found on numerical study of this model but it may be noted that the FEM convergence analysis, using this simplest possible approach, has not been reported anywhere of this modified Navier Stokes equations coupled with energy equation in the context of HNF flow. For finite element analysis we consider a domain $\Omega \subset {R}^2$ with its boundary $\partial \Omega$ which is Lipschitz continuous. The existence uniqueness of weak solution of a similar type of model on a Lipschitz continuous boundary was proved in the paper\cite{Yyh11,CZ95}. In the above literature mostly we have observed the use of one particular type of nanofluid. But here we have used the hybrid nanofluid $(Cu, Al_2O_3)$/water to explore more on heat conduction and heat flow and enhancement of heat transfer in this case  than that of the situation in case of a  nano fluid. Next to that we focus on exploring the influence of this  Hybrid nano particle on the fluid flow in intricate domains like L-shaped and H-shaped domains with partial/total heating of left vertical wall. Even finding the numerical solution of this problem is a difficult task. The presence of temperature field in the equations made it more challenging to reach the convergence of the numerical scheme for a desired range of parametric set up. Fortunately, for an intended range of parameters we have convergence of the numerical scheme used here which is discussed in the Numerical Experiment part. Numerical scheme is validated by making a comparison with existing result and the output data is presented through Matlab plotting and tabulated data.\\ 

The arrangement of the paper is as follows. Section 3.1 and section 3.2 deal with the  non-dimensionalized and vector form of the model respectively. Section 4 contains a detailed description of weak formulation of the model whereas section 5 contains the discrete formulation of the problem. In section 6 we derive the apriori finite element error estimate to establish the convergence of the finite element scheme and there after go on applying the scheme to analyse the natural convection in complex domains. In section 7 we have focused on Numerical Experiment. In section $7.1$  derivation of elementary matrix formulation is presented. In $7.2$ in table 1 establishes the validity of our numerical scheme as the present results match nicely with those in the existing literature. Also in section 7.3 we have carried out grid validation  test on both H and L shape domains. We have presented a few plots obtained from numerical experiments done on H-shape domain for  the clear fluid along with the output data in table 4 to validate our numerical scheme in section 8.1. Those plots depict patterns as reported in the literature \cite{LSH03}. We have presented the obtained results by tabulated values and by a  few selected plots to explain the behaviour of flow pattern. The effect of pertinent parameters like Prandtl number $(1\leq Pr\leq 10)$, volume fraction $( 0\leq \phi\leq 0.01)$, and Rayleigh number $(1\leq Ra\leq 10^{5})$ are described also by Matlab plotting and tabular representations in section 8.2 and 8.3. In section 9 numerical experiment is carried over an L shaped domain, a particular case of $\Omega$, where we have explained how adding of nano particles enhances heat conduction and a very clear comparison of Nusselt number between clear fluid(CF) and HNF for different type of volume fraction$(0.001\leq \phi \leq 0.01)$ is depicted by Matlab plotting and by tabulated values. Also the effect of the changes of input parameters like Rayleigh number$(Ra)$, and source length respectively can be observed by the corresponding changes in streamlines and isotherms plots. At the end we have tried to explain the  direct effect of the length of heat source on the Nusselt number which is presented in tabular form in section 9.3. The major inferences are concluded in section 10.
Therefore, the above stated numerical experiments give us an opportunity to show that adding of nano particles really enhances the low conduction flow problem as desired by industries. 
\section{Nomenclature}\label{sec:1}

x,y  Cartesian coordinates $(m)$\vspace{2mm}.\\
Pr  Prandtl number.\vspace{2mm}\\
Ra Rayleigh number.\vspace{2mm}\\
T  Temperature,$(K)$.\vspace{2mm}\\
$\theta$ Non-dimensional temperature.\vspace{2mm}\\
U dimensionless velocity component along x direction.\vspace{2mm}\\
V dimensionless velocity component along y direction.\vspace{2mm}\\
$\alpha$ thermal diffusivity, $(m^2/s)$.\vspace{2mm}\\
$\beta$ thermal expansion coefficient,$(K^{-1})$.\vspace{2mm}\\
$\mu$ dynamic viscosity,$(Pa.s)$.\vspace{2mm}\\
$\rho$ density,$(Kg/m^3)$.\vspace{2mm}\\
\texttt{Abbreviations:}\vspace{2mm}\\
CF $\rightarrow$ fluid \vspace{2mm}\\
V.F. $\rightarrow$ volume fraction \vspace{2mm}\\
HNF $\rightarrow$ hybridnano fluid.\vspace{4.4mm}\\
\section{About Our Model}
\label{}
Let $\Omega \subset R^2$ be a bounded and connected domain with Lipschitz continuous boundary $\partial{\Omega}.$ We assume that the HNF used here is non Newtonian, homogeneous, laminar and incompressible in nature. The governing equations consists of the continuity equations, momentum equations in both x and y directions and energy equation as written below.\vspace{3mm}
\subsection{Non-dimensionalized form \cite{CSM19}:} 
\begin{align*}
\frac{\partial U}{\partial x} + \frac{\partial V}{\partial y}=0\\ 
-(\frac{\rho_f}{\rho_{hnf}})(\frac{\mu_{hnf}}{\mu_f})Pr(\frac{\partial^2 U}{\partial x^2}+\frac{\partial^2 V}{\partial y^2} ) +(U \frac{\partial U}{\partial x}+V \frac{\partial U}{\partial y}) +\frac{\rho_f}{\rho_{hnf}}\frac{\partial p}{\partial x}=0\\
-(\frac{\rho_f}{\rho_{hnf}})(\frac{\mu_{hnf}}{\mu_f})Pr(\frac{\partial^2 V}{\partial x^2}+\frac{\partial^2 V}{\partial y^2} ) +(U \frac{\partial V}{\partial x}+V \frac{\partial V}{\partial y}) +\frac{\rho_f}{\rho_{hnf}}\frac{\partial p}{\partial y}-\frac{(\rho \beta)_{hnf}}{\rho_{hnf}\beta_f}Pr Ra T=0 \\
U\frac{\partial T}{\partial x}+V \frac{\partial T}{\partial y} -\frac{\alpha_{hnf}}{\alpha_f} (\frac{\partial^2 T}{\partial x^2}+\frac{\partial^2 T}{\partial y^2})=0, 
\end{align*}
${U}={0},~V=0$ on $\partial {\Omega}$, 
 $T=1$ on a part of $\partial {\Omega}$ and $T=0$ elsewhere on $\partial {\Omega}$.

\subsection{The vector form of the above problem:} 
Let, $\textbf{u}=(U,V)$.
\begin{align}
div ~\textbf{u}=0\\ 
-(\frac{\rho_f}{\rho_{hnf}})(\frac{\mu_{hnf}}{\mu_f})Pr\triangle \textbf{u} + (\textbf{u} \cdot \nabla)\textbf{u}+(\frac{\rho_f}{\rho_{hnf}})\nabla p= \frac{(\rho \beta)_{hnf}}{\rho_{hnf}\beta_f}Pr Ra jT~ \\,j=(0,1)\nonumber \\
-\frac{\alpha_{hnf}}{\alpha_f}\triangle T +\textbf{u} \cdot \nabla T =0
 \end{align}
$\textbf{u}=\textbf{0}$; 
$T=1$ on a part of $\partial {\Omega}$ and $T=0$ elsewhere on $\partial {\Omega}$.
\section{The Weak Formulation:}
\textbf{Definition:} We need the following spaces to define the weak formulation as well as to establish the finite element error estimation\vspace{2.4mm}\\$H^1(\Omega)=\{v~\in L^2(\Omega):\frac{\partial v}{\partial x},\frac{\partial v}{\partial y}~\in~L^2(\Omega)\}$,\vspace{2.4mm}\\
$V=\{\texttt{v}~\in~H^1(\Omega)^2:\texttt{v}=\textbf{0}~\mbox{on}~\partial\Omega\}$ .\vspace{2.4mm}\\ $Z=H^1(\Omega), \hspace{2.3mm} W=H_{0}^1(\Omega)=\{v~\in~H^1(\Omega):v=0~\mbox{on}~\partial\Omega\}.$\vspace{2.4mm}\\The norm on $H_0^1(\Omega)$ is given by $\|\nabla v\|_2=\mid v\mid_{1,\Omega}$.\vspace{2.6mm}\\
Let us introduce the space $X=\{v\in V :\mbox{div}~v=0\}$.\vspace{2.4mm}\\
$$Q=L_{0}^2(\Omega)=\{v\in L^2(\Omega):\int_{\Omega} v=0\},$$\\
We introduce the following semi-norm on the Space $H^{m}(\Omega)$ as follows
$$\mid u\mid_{m}=(\Sigma_{\mid \alpha \mid=m}\int_{\Omega}\mid \partial ^{\alpha}u(x)\mid^{p} dx)^{1/p}$$.\vspace{3mm}\\
\textbf{Note 1:} Notation $(.,.)$  stands for the inner product in $L^2(\Omega)$ or $(L^2(\Omega))^2$.\vspace{3mm}\\
\textbf{Note 2:} We will prove our theorems with the boundary conditions $ \textbf{u}=\textbf{0},T=0~ on~ \partial{\Omega}$ and adding a source term $f\in L^2(\Omega)$. As the case with non zero temperature condition on a part of the boundary is an easy consequence of this by using translation.\vspace{2.4mm}\\
Let us apply Green's formula to obtain the variational formulation\\ of the above  problem and rename the problem as problem  $(P)$.\vspace{3mm}\\
To find $(\textbf{u},p,T)~\in~V \times Q \times W$ such that
\begin{align}
a_0 (\textbf{u},\textbf{v}) +a_1 (\textbf{u};\textbf{u},\textbf{v})-b(\textbf{v},p)- a_2 (T,\textbf{v})+a_3 (T,\theta)
\nonumber\\ 
+a_4(\textbf{u};T,\theta)=(f,\theta);~~b(\textbf{u},q)=
0.
\end{align}
 $~\forall~(\texttt{v},q,\theta)~\in~V \times Q \times W$, $f \in L^2(\Omega)$ where\\
 \begin{align*}
 a_0(\textbf{u},\textbf{v})=(\frac{\rho_f}{\rho_{hnf}})(\frac{\mu_{hnf}}{\mu_f})Pr\int_{\Omega}\nabla \textbf{u} \cdot \nabla \textbf{v}\\
a_2 (T,\textbf{v})=\frac{(\rho \beta)_{hnf}}{\rho_{hnf}\beta_f}Pr Ra \int_{\Omega}(j T ) \cdot \textbf{v}\\
a_3(T,\theta)= \frac{\alpha_{hnf}}{\alpha_f} \int_{\Omega} \nabla T \cdot \nabla \theta\\
b(\textbf{u},p)=\int_{\Omega}(\frac{\rho_f}{\rho_{hnf}})(\nabla \cdot \textbf{u} )p
\\
a_1 (\textbf{u};\textbf{u},\textbf{v})= \int_{\Omega} (\textbf{u} \cdot \nabla \textbf{u})\cdot \textbf{v}\\
 a_4(\textbf{u};T,\theta)=\int_{\Omega} (\textbf{u} \cdot \nabla T)\theta\\
(f,\theta)=\int_{\Omega} f \theta
\end{align*}\vspace{2.3mm}\\
 
In order to estimate the first nonlinear term  $a_1 (\textbf{u};\textbf{u},\textbf{v})$ and $a_4 (\textbf{u};T,\theta)$, we introduce the following tri-linear forms on $(H_0^1(\Omega)^2)^3$ and $(H_0^1(\Omega))^3$ respectively as follows:\\ 
\begin{equation}
a_1(\textbf{w};\textbf{u},\textbf{v})
= \sum_{i,j=1}^2\int_\Omega (\textbf{w}_j(\frac{\partial \textbf{u}_i}{\partial x_j})\textbf{v}_i) dx~~
\end{equation}\vspace {3mm}\\
\begin{equation}
\begin {split}
a_4(\textbf{u};T,\theta)
&= \sum_{j=1}^2 \int_\Omega (\textbf{u}_j(\frac{\partial T}{\partial x_j})\theta) dx
\end{split}
\end{equation} \vspace {3mm}\\
Our tri-linear terms satisfy the following continuity condition on 
\begin{equation}
\begin {split}
\mid a_1(\textbf{w};\textbf{u},\textbf{v}) \mid \leq N_1 {{\mid \textbf{u} \mid}_{1,\Omega} {\mid \textbf{v} \mid}_{1,\Omega} {\mid \textbf{w} \mid}_{1,\Omega}} &\\
\mid a_2(\textbf{u};T,\theta) \mid \leq N_2{{\mid \textbf{u} \mid}_{1,\Omega} {\mid T \mid}_{1,\Omega} {\mid \theta \mid}_{1,\Omega}} 
  \end{split}
  \end{equation} 
\vspace {1.1mm}\\,where $N_1,N_2\geq 0$ depend on $\Omega$\cite{SGS10}.
\vspace{4.2mm}\\
\texttt{Inf-Sup Condition:}
Our velocity and pressure space $(H_0^1(\Omega))^2 \mbox{and}~ L_0^2(\Omega)$ satisfies the well known Babuska-Brezzi condition \cite{Kes88}  
$\exists$ a constant $A> 0$ such that\vspace{1.6mm}\\
\begin{equation}
\underset{q \in Q}{inf}\underset{\textbf{v} \in V}{ sup}\frac{(\nabla \cdot \textbf{v},q)}{{\mid \textbf{v} \mid}_{1,\Omega} {\| q \|_2}} \geq {A} > 0.
\end{equation}\vspace{1.3mm}\\
which is equivalent to \cite{GVR09}\vspace{1.3mm}\\
\begin{equation}
\underset{q \in Q}{inf}\underset{(\textbf{u},\theta)\in V \times W}{ sup}\frac{(\nabla \cdot \textbf{v},q)}{\|(\textbf{v},\theta)\|_{V \times W} \|q\|_2} \geq BA >0,\end{equation}
$B>0$ is a constant.\vspace{4.2mm}\\
\texttt{Remark:}The existence uniqueness and stability bound of the considered model can be derived by following the same proof presented in \cite{CZ95,Yyh11}.\vspace{2.1mm}\\
\section{Finite element formulation:}\label{sec:4}
\label{}
\subsection{Finite Element Spaces:} Let $\tau_h$ be the finite element partition of the domain $\Omega$ with mesh size h as $\Omega=\underset{k \in \tau_h}{\cup K }$. Define the finite element spaces corresponding to the Velocity, pressure and temperature respectively as\vspace{3.2mm}\\ 
$X_h=\{\textbf{u}_h~\in V\cap (C^{0}(\Omega))^{2}:V_h\mid_{K} \in (P_{2}(K))^{2},\forall K \in \tau_h \}$\vspace{2.3mm}\\ 
$Q_h=\{q_h \in Q:{q_{h}}\mid_{K}~\in~P_{1}(K),\forall K~\in \tau_h\}\vspace{2.3mm}\\
R_h=\{T_h \in W \cap C^{0}(\Omega) : T_h\mid_{K}~\in P_2(K)~\forall K \in \tau_h)\}$\vspace{2.3mm}\\
$V_h=\{ \textbf{v}_h\in X_h:\int_{\Omega} ( q_h \nabla) \cdot \textbf{v}_h=0,\forall ~q_h \in Q_h\}$\vspace{2.1mm}\\
The discrete form of the inf-sup condition is as follows
\begin{equation}
\exists ~a~constant~L~> 0~s.t.~\underset{q_h \in Q_h}{inf}\underset{\textbf{v}_h \in X_h}{ sup}\frac{(\nabla \cdot \textbf{v}_h,q_h)}{{\mid \textbf{v}_h \mid}_{1,\Omega} {\| q_h \|_2}} \geq {L} > 0
\end{equation}
where $P_r(K)$ is the space of polynomial of degree $r$.\vspace{4mm}\\
We assume that the discrete analogue of equivalence of $(1)$ and $(13)$ also holds
\begin{equation}
\underset{q_h \in Q_h}{inf}\underset{(\textbf{v}_h,\theta_h)\in X_h \times R_h}{ sup}\frac{(\nabla \cdot \textbf{v}_h,q_h)}{\|(\textbf{v}_h,\theta_h)\|_{X_h \times R_h} \|q_h\|_2} \geq BL >0,\end{equation}
$B>0$ is a constant.\vspace{4mm}\\

As the finite element spaces are conforming, our weak form $(8)$ also satisfies the following equations. \vspace{2.2mm}\\
$\forall ~\textbf{v}_h~\in ~X_h,\theta_h ~\in~R_h,q_h~\in~Q_h$ \\
\begin{equation}
a_0 (\textbf{u},\textbf{v}_h) +a_1 (\textbf{u};\textbf{u},\textbf{v}_h)-b(\textbf{v}_h,p)- a_2 (T,\textbf{v}_h)+a_3 (T,\theta_h)\\ \nonumber \end{equation}
\begin{equation}
+a_4(\textbf{u};T,\theta_h)=(f,\theta_h);~~b(\textbf{u},q_h)=0.
\end{equation}\vspace{2mm}\\
\textbf{Discrete Problem:} We consider the weak formulation of the discrete problem as :~
to find $(\textbf{u}_h,p_h.T_h)\in X_h \times Q_h \times R_h$
\begin{equation}
a_0 (\textbf{u}_h,\textbf{v}_h) +a_1 (\textbf{u}_h;\textbf{u}_h,\textbf{v}_h)-b(\textbf{v}_h,p)- a_2 (T_h,\textbf{v}_h)+a_3 (T_h,\theta_h)\\ \nonumber \end{equation}
\begin{equation}
+a_4(\textbf{u}_h;T_h,\theta_h)=(f,\theta_h);~~b(\textbf{u}_h,q_h)=0.
\end{equation}\vspace{1mm}\\$\forall~(\textbf{v}_h,q_h,\theta_h)~\in~X_h \times Q_h \times R_h$.\vspace{4mm}\\As the finite element spaces are conforming we have existence uniqueness of weak solution of the discrete problem under the same conditions and the solution is stable.\\ 
\section{Error estimation:}
\texttt{Interpolation error estimate \cite{SGS10}:} We assume that our weak solution $(\textbf{u},p,T)\in (X \cap (H^{3}(\Omega))^{2}) \times (Q \cap H^{2}(\Omega)) \times (W \cap H^{3}(\Omega))$. Then for the finite element spaces $(X_h,Q_h,R_h)$ we have the following interpolation error estimate\vspace{2.1mm}\\
\begin{align}
\underset{\textbf{v}_h \in X_h}{inf}\{ \|\textbf{u}-\textbf{u}_h \|_{0}+h\|\nabla(\textbf{u}-\textbf{v}_h)\|_{0}\}\leq Ch^{3}\mid \textbf{v} \mid_{3}\vspace{2.1mm}\\
\underset{T_h \in R_h}{inf}\{ \|T-T_h \|_{0}+h\|\nabla(T-T_h)\|_{0}\}\leq Ch^{3}\mid T\mid_{3}\vspace{2.1mm}\\
\underset{q_h \in Q_h}{inf}\{ \|p-q_h \|_{0}\}\leq Ch^{2}\mid p\mid_{2}.
\end{align}
\begin{theorem}
Let $(\textbf{u},p,T)$ be the unique solution of the modified Navier-Stokes equations coupled with energy equation model in the context of hybrid nano-fluid. 
Let the discrete problem has the weak solution $(\textbf{u}_h,p_h,T_h)$. Let R be the stability bound of the weak solutions which depend on the norm of source term. \vspace{2mm}\\
Then,\vspace{6mm}\\ $${\mid \textbf{u}-\textbf{u}_h \mid}_{1,\Omega} \leq \frac{C_1 h^2}{\surd{L_1}}(\surd{K_1} {\mid \textbf{u} \mid}_{3}+\surd{K_3}{\mid T \mid}_{3})+C_2h^{2}\surd{K_2}{\mid p \mid}_{2}),$$ \vspace{2mm}\\
$${\mid T-T_h \mid}_{1,\Omega} \leq \frac{C_1 h^2}{\surd{L_2}}(\surd{K_1}{\mid \textbf{u} \mid}_{3}+\surd{K_3}{\mid T \mid}_{3})+C_2h^{2}\surd{K_2}{\mid p \mid}_{2})$$\vspace{3mm}\\ where $L_1=\big\{\frac{(\frac{\rho_f}{\rho_{hnf}})(\frac{\mu_{hnf}}{\mu_f})Pr}{2}-(N_1+N_2)R-{\frac{(\rho \beta)_{hnf}}{\rho_{hnf}\beta_f} Pr Ra c^2}\},\\L_2=\big\{\frac{{\alpha}_{hnf}}{2\alpha_f}-\frac{(\rho \beta)_{hnf}}{2\rho_{hnf}\beta_f} Pr Ra c^2-5N_2R\},$ provided $ L_1,L_2$ are strictly positive.\vspace{2mm}\\
 Here, $$K_1=(\frac{3(\frac{\rho_f}{\rho_{hnf}})(\frac{\mu_{hnf}}{\mu_f})Pr}{2} +\frac{6(N_1R)^2}{(\frac{\rho_f}{\rho_{hnf}})(\frac{\mu_{hnf}}{\mu_f})Pr}+N_2R),$$ \\
$$ K_2=\frac{\frac{3}{2}(\frac{\rho_f}{\rho_{hnf}})}{(\frac{\mu_{hnf}}{\mu_f}) Pr}, K_3=(N_2R+\frac{\frac{(\rho \beta)_{hnf}}{\rho_{hnf}\beta_f}Pr Ra c^2}{2}+\frac{\alpha}{2}).$$
 \end{theorem}
Proof: We start with the following difference equation by subtracting $(13)$ from equation $(12)$ 
\begin{equation}
a_0 (\textbf{u}-\textbf{u}_{h},\textbf{v}_h) +a_1 (\textbf{u};\textbf{u},\textbf{v}_h)- a_1 (\textbf{u}_{h};\textbf{u}_{h},\textbf{v}_h)-b(\textbf{v}_h,p-p_h)- a_2 (T-T_{h},\textbf{v}_h)\\ \nonumber \end{equation}
\begin{equation}
+a_3 (T-{T}_h,\theta_h)+a_4(\textbf{u};T,\theta_h)-a_4(\textbf{u}_{h};T_h,\theta_h)=0,~~b(\textbf{u}-\textbf{u}_h,q_h)=0
\end{equation}
$\forall ~\textbf{v}_h~\in ~X_h,\theta_h ~\in~R_h,q_h~\in~Q_h$ \vspace{3mm}\\

Let $\boldsymbol{\xi}_h =\textbf{u}_h- {\boldsymbol{\phi}}_h,~~~~~\eta_h=p_h-P_h ,~~~~~\delta_h=T_h-\tau_h.$\\ 
Denote by $\boldsymbol{\xi} = \textbf{u}- \boldsymbol{\phi}_h,~~~~\eta = p-P_h ,~~~~~\delta=T- \tau_h$\\ for any function $\boldsymbol{\phi}_h\in V_h , ~P_h \in Q_h, ~\tau_h \in R_h$ .\vspace{3mm}\\

Then taking $\boldsymbol{v}_h=\boldsymbol{\xi}_h\in V_h,~q_h=\eta_h\in Q_h,~\mbox{and}~\theta_h=\delta_h\in R_h$ after using the above introduced notations in the above difference equation $(17)$ we have\\
\begin{equation}
a_0(\boldsymbol{\xi}_h,\boldsymbol{\xi}_h)+a_3(\delta_h,\delta_h)=a_2(\delta_h,\boldsymbol{\xi}_h)+a_1(\textbf{u};\textbf{u},\boldsymbol{\xi}_h)-a_1(\textbf{u}_h;\textbf{u}_h,\boldsymbol{\xi}_h)+a_4(\textbf{u};T,\delta_h)\nonumber \end{equation}
\begin{equation}
-a_4(\textbf{u}_h;T_h,\delta_h)+a_0(\boldsymbol{\xi},\boldsymbol{\xi}_h)-b(\boldsymbol{\xi}_h,\eta)+a_2(\delta_h,\boldsymbol{\xi}_h)
-a_2(\delta,\boldsymbol{\xi}_h)+a_3(\delta,\delta_h),
\nonumber\end{equation}
\begin{equation}
b(\boldsymbol{\xi}_h,\eta_h)=b(\boldsymbol{\xi},\eta_h)=0
\end{equation}\vspace{2.3mm}\\
Let us evaluate the difference of the first non-linear term as
\begin{equation}
a_1(\textbf{u};\textbf{u},\boldsymbol{\xi}_h)-a_1(\textbf{u}_h;\textbf{u}_h,\boldsymbol{\xi}_h)
=a_1(\textbf{u};\boldsymbol{\xi},\boldsymbol{\xi}_h)+a_1(\boldsymbol{\xi};\textbf{u}_h,\boldsymbol{\xi}_h)-a_1(\boldsymbol{\xi}_h;\textbf{u}_h,\boldsymbol{\xi}_h)
\nonumber \end{equation}
Therefore,
\begin{equation}
\begin{split}
&\mid a_1(\textbf{u};\textbf{u},\boldsymbol{\xi}_h)-a_1(\textbf{u}_h;\textbf{u}_h,\boldsymbol{\xi}_h)\mid\\
&\leq N_1(\|\nabla \textbf{u}\|_2+\|\nabla \textbf{u}_h\|_2)\|\nabla \boldsymbol{\xi}\|_2\|\nabla \boldsymbol{\xi}_h\|_2
+N_1\|\nabla \textbf{u}_h\|_2 \|\nabla \boldsymbol{\xi}_h\|_2^2\\
&\leq 6\frac{{N_1}^2 R^{2}}{(\frac{\rho_f}{\rho_{hnf}})(\frac{\mu_{hnf}}{\mu_f})Pr}\|\nabla \boldsymbol{\xi}\|_2^{2}+\frac{((\frac{\rho_f}{\rho_{hnf}})(\frac{\mu_{hnf}}{\mu_f})Pr}{6}\|\nabla \boldsymbol{\xi}_h\|_2^{2}+N_1R\|\nabla \boldsymbol{\xi}_h\|_2^{2}
\end{split}
\end{equation}
Let us evaluate the second nonlinear term as follows
\begin{equation}
\begin{split} 
&\mid a_4(\textbf{u};T,\delta_h)-a_4(\textbf{u}_{h};T_h,\delta_h)\mid\\ &\leq\int_{\Omega}\mid[(\textbf{u}-\textbf{u}_h) \cdot \nabla T - \textbf{u}_h \cdot \nabla (T-T_h)]\delta_h\mid\\ 
&\leq \int_{\Omega} \mid [(\boldsymbol{\xi}-\boldsymbol{\xi}_h) \cdot \nabla T]\delta_h \mid +\int_{\Omega}\mid [u_h \cdot \nabla (\delta-\delta_h)]\delta_h \mid
\end{split}
\end{equation}\vspace{3mm}\\
Let us introduce  $\boldsymbol{\xi}=({\xi}^1,{\xi}^2),\boldsymbol{\xi}_h=({\xi}_h^1,{\xi}_h^2)$\vspace{2.3mm}\\ So, $\|\nabla \boldsymbol{\xi}\|_2^2=\|\nabla {\xi}^1\|_2^2+\|\nabla{\xi}^2\|_2^2$.\vspace{2.3mm}\\
$\|\nabla \boldsymbol{\xi}_h\|_2^2=\|\nabla {\xi}_h^1\|_2^2+\|\nabla {\xi}_h^2\|_2^2$\vspace{2.6mm}\\

The first term in equation $(20)$ can be evaluated as\\
\begin{equation}
\begin{split}
& \int_{\Omega}\mid [(\boldsymbol{\xi}-\boldsymbol{\xi}_h) \cdot \nabla T]\delta_h \mid \\
&\leq \int_{\Omega}\mid ({\xi}^1-{\xi}_h^1)\frac{\partial T}{\partial x}\delta _h \mid +\int_{\Omega}\mid ({\xi}^2-{\xi}_h^2)\frac{\partial T}{\partial y}\delta _h \mid \\
&\leq N_2\|\nabla T\|_2\big\{(\| \nabla {\xi}^1\|_2+\|\nabla { \xi}_h^1\|_2)\|\nabla \delta_h\|_2+(\| \nabla {\xi}^2\|_2+\|\nabla {\xi}_h^2\|_2)\|\nabla \delta_h\|_2\}\vspace{2.3mm}\\
&\leq N_2 R(\|\nabla \boldsymbol{ \xi}\|_2^2+\|\nabla \boldsymbol{\xi}_h\|_2^2+2\|\nabla \delta_h\|_2^2)\end{split}
 \end{equation}
Similarly, the second term in equation $(20)$ can be estimated as\\
\begin{equation}
 \int_{\Omega}\mid[(\textbf{u}_h \cdot \nabla (\delta-\delta_h)]\delta_h \mid 
\leq N_2R[\|\nabla \delta\|_2^2+3\|\nabla \delta_h\|_2^2]
\end{equation}\vspace{4mm}\\

The estimation of the linear terms are as follows\\
\begin{equation}
\begin{split}
\mid a_0(\boldsymbol{\xi},\boldsymbol{\xi}_h)\mid &\leq (\frac{\rho_f}{\rho_{hnf}})(\frac{\mu_{hnf}}{\mu_f})Pr\int_{\Omega}\mid \nabla \boldsymbol{ \xi} \cdot \nabla \boldsymbol{\xi}_h\mid\\
&\leq \frac{((\frac{\rho_f}{\rho_{hnf}})(\frac{\mu_{hnf}}{\mu_f})Pr }{2}(\frac{1}{3}\|\nabla \boldsymbol{\xi}_h \|_2^2+3\|\nabla \boldsymbol{\xi} \|_2^2)
\end{split}    
\end{equation}
\vspace{3mm}\\
\begin{equation}
\begin{split}
\mid b(\boldsymbol{\xi}_h,\eta)\mid &\leq \int_{\Omega}(\frac{\rho_f}{\rho_{hnf}})
\mid (\nabla \cdot \boldsymbol{\xi}_h )\eta\mid\\
&\leq (\frac{\rho_f}{\rho_{hnf}})\|\nabla \cdot \boldsymbol{\xi}_h \|_2 \|\eta\|_2\\
&\leq (\frac{\rho_f}{2\rho_{hnf}})((\frac{\mu_{hnf}}{3\mu_f}Pr)\|\nabla \boldsymbol{\xi}_h\|_2^2+\frac{3}{(\frac{\mu_{hnf}}{\mu_f})Pr}\|\eta \|_2^2)
\end{split}
\end{equation}
\vspace{3mm}\\
\begin{equation}
\begin{split}
\mid a_2 (\delta_h,\boldsymbol{\xi}_h)\mid &\leq \frac{(\rho \beta)_{hnf}}{\rho_{hnf}\beta_f}Pr Ra \int_{\Omega}\mid (j \delta_h ) \cdot \boldsymbol{\xi}_h \mid \\
&\leq \frac{(\rho \beta)_{hnf}}{\rho_{hnf}\beta_f}Pr Ra c^2\|\nabla \delta_h \|_2 \|\nabla \boldsymbol{\xi}_h\|_2\\
&\leq \frac{\frac{(\rho \beta)_{hnf}}{\rho_{hnf}\beta_f}Pr Ra c^2}{2} (\|\nabla \delta_h\|_2^2+\|\nabla \boldsymbol{\xi}_h\|_2^2)
\end{split}
\end{equation}\vspace{3mm}\\
\begin{equation}
\begin{split}
\mid a_2 (\delta,\boldsymbol{\xi}_h)\mid &\leq \frac{(\rho \beta)_{hnf}}{\rho_{hnf}\beta_f}Pr Ra \int_{\Omega}\mid (j \delta ) \cdot \boldsymbol{\xi}_h \mid \\
&\leq \frac{\frac{(\rho \beta)_{hnf}}{\rho_{hnf}\beta_f}Pr Ra c^2}{2} (\|\nabla \delta\|_2^2+\|\nabla \boldsymbol{\xi}_h\|_2^2)
\end{split}
\end{equation}\vspace{3mm}\\
\begin{equation}
\begin{split}
\mid a_3(\delta,\delta_h)\mid &\leq \int_{\Omega} \mid \frac{\alpha_{hnf}}{\alpha_f} \nabla \delta  \cdot \nabla \delta_h \mid \\
&\leq \frac{\frac{\alpha_{hnf}}{\alpha_f}}{2}(\|\nabla \delta_h\|_2^2+\|\nabla \delta\|_2^2)
\end{split}
\end{equation}
[{Using the stability estimate, Young's inequality, continuity of the 2nd non-linear term.}]\\
Using the estimates (19)-(27) in equation (18) we have
\begin{equation}
L_1 \| \nabla \boldsymbol{\xi}_h \|_2^2
+L_2\| \nabla \delta_h\|_2^2
\leq ({K_1}\|\nabla \boldsymbol{\xi}\|_2^2+{K_2}\|\eta \|_2^2+{K_3}\|\nabla \delta\|_2^2 )
\end{equation}

where
$$K_1=(\frac{(\frac{3\rho_f}{\rho_{hnf}})(\frac{\mu_{hnf}}{\mu_f})Pr}{2} +\frac{6(N_1R)^2}{(\frac{\rho_f}{\rho_{hnf}})(\frac{\mu_{hnf}}{\mu_f})Pr}+N_2R),$$ \\
$$ K_2=\frac{\frac{3}{2}(\frac{\rho_f}{\rho_{hnf}})}{(\frac{\mu_{hnf}}{\mu_f}) Pr}, K_3=(N_2R+\frac{\frac{(\rho \beta)_{hnf}}{\rho_{hnf}\beta_f}Pr Ra c^2}{2}+\frac{\alpha}{2}).$$\vspace{4mm}
 \begin{equation}
\Rightarrow \| \nabla \boldsymbol{\xi}_h \|_2 \leq \frac{1}{2\surd{L_1}}(\surd{K_1}\|\nabla \boldsymbol{\xi}\|_2+\surd{K_2}\|\eta \|_2+\surd{K_3}\|\nabla \delta\|_2) \end{equation}
 \begin{equation}\mbox{and}~\| \nabla \delta_h \|_2 \leq \frac{1}{2\surd{L_2}}(\surd{K_1}\|\nabla \boldsymbol{\xi}\|_2+\surd{K_2}\|\eta \|_2+\surd{K_3}\|\nabla \delta\|_2).
\end{equation}\\
 where $L_1=\big\{\frac{(\frac{\rho_f}{\rho_{hnf}})(\frac{\mu_{hnf}}{\mu_f})Pr}{2}-(N_1+N_2)R-{\frac{(\rho \beta)_{hnf}}{\rho_{hnf}\beta_f} Pr Ra c^2}\},\\L_2=\big\{\frac{{\alpha}_{hnf}}{2\alpha_f}-\frac{(\rho \beta)_{hnf}}{2\rho_{hnf}\beta_f} Pr Ra  c^2-5N_2R\},$
 provided $L_1,L_2>0$\\ i.e. are strictly positive.\vspace{5mm}\\
Therefore, using interpolation error estimate in (29) and (30)  we have
 \begin{equation}
 \begin{split}
& \|\nabla(\textbf{u}-\textbf{u}_h)\|_2\\
 &\leq \|\nabla \boldsymbol{\xi}\|_2+\|\nabla \boldsymbol{ \xi}_h\|_2\\
 & \leq  \|\nabla \boldsymbol{\xi} \|_2+\frac{1}{2\surd{L_1}}(\surd{K_1}\|\nabla \boldsymbol{ \xi}\|_2+\surd{K_2}\|\eta \|_2+\surd{K_3}\|\nabla \delta\|_2)\\
&= (1+\frac{\surd{K_1}}{2\surd{L_1}})\|\nabla(\textbf{u}-\boldsymbol{\phi}_h)\|_2+\frac{\surd{K_2}}{2\surd{L_1}}\|(p-P_h )\|_2+\frac{\surd{K_3}}{2\surd{L_1}}\|\nabla(T-\tau_h)\|_2\\
 &\leq c_1 h^2 \mid \textbf{u} \mid_{3}+c_2 h^2 \mid p \mid_{2}+c_3 h^2 \mid T\mid_{3}
 \end{split}
 \end{equation}\vspace{3.6mm}\\
 \begin{equation}
 \begin{split}
& \|\nabla(T-T_h)\|_2\\
 &\leq \|\nabla \delta\|_2+\|\nabla \delta_h\|_2\\
 & \leq  \|\nabla \delta\|_2+\frac{1}{2\surd{L_1}}(\surd{K_1}\|\nabla \boldsymbol{\xi}\|_2+\surd{K_2}\|\eta \|_2+\surd{K_3}\|\nabla \delta\|_2)\\
&=(1+\frac{\surd{K_1}}{2\surd{L_1}})\|\nabla(
\textbf{u}- \boldsymbol{\phi}_h)\|_2+\frac{\surd{K_2}}{2\surd{L_1}}\|(p-P_h )\|_2+\frac{\surd{K_3}}{2\surd{L_1}}\|\nabla(T-\tau_h)\|_2\\
 &\leq c_1 h^2 \mid \textbf{u} \mid_{3}+c_2 h^2 \mid p \mid_{2}+c_3 h^2 \mid T \mid_{3}
 \end{split}
 \end{equation}
\vspace{4mm}\\
Now we concentrate on the pressure approximation.
\begin{theorem}
Let $(\textbf{u},p,T)$ be the unique solution of the modified Navier-Stokes equations coupled with energy equation model in the context of hybrid nano-fluid. 
Let the discrete problem has the weak solution $(\textbf{u}_h,p_h,T_h)$. Let R be the stability bound of the weak solutions which depend on the norm of source term.\vspace{3mm}\\ Then,
\begin{equation}
 \begin{split}
 &\|p-P_h\|_2\\
 &\leq \frac{\rho_{hnf}}{BL\rho_{f}}[((\frac{\rho_f}{\rho_{hnf}})(\frac{\mu_{hnf}}{\mu_f})Pr+2N_1R+N_2R) 
 +(\frac{(\rho \beta)_{hnf}}{ \rho_{hnf}\beta_f}Pr Ra c^2  +N_2R+\frac{\alpha_{hnf}}{\alpha_f})]\\
 &\times(c_1 h^2 {\mid \textbf{u} \mid}_{3}+c_2 h^2 {\mid p\mid}_{2}+c_3 h^2 {\mid T \mid}_{3}) +(1+(BL)^{-1})Ch^2\mid p \mid _2.
\end{split}
 \end{equation}\end {theorem}
\textbf{Proof:}
For pressure approximation we use the equivalence of inf-sup condition in discrete form. Here using $(11)$ we have\vspace{4.5mm}\\
\begin{equation}
\begin{split}
&BL\|p_h-P_h\|_2\\ &\leq \underset{(\textbf{v}_h,\theta_h) \in X_h \times R_h}{sup} \frac{(\nabla \cdot \textbf{v}_h,p_h-P_h)}{\|(\textbf{v}_h,\theta_h)\|}\\
&= \underset{(\textbf{v}_h,\theta_h) \in X_h \times R_h}{sup} \frac{b(\textbf{v}_h,p_h-p)+b(\textbf{v}_h,p-P_h)}{\| (\textbf{v}_h,\theta_h )\|}
\end{split}
\nonumber \end{equation}
\begin{equation}
\leq \underset{\textbf{v}_h \in X_h}{sup} \frac{\frac{\rho_{hnf}}{\rho_{f}}\mid[a_0(\textbf{u}-\textbf{u}_h,\textbf{u}_h)+a_1(\textbf{u};\textbf{u},\textbf{u}_h)-a_1(\textbf{u}_h;
\textbf{u}_h,\textbf{v}_h)-a_2(j(T-T_h),\textbf{v}_h)]\mid}{\mid \textbf{v}_h \mid_{1,\Omega}}
\nonumber\end{equation}
\begin{equation}
+ \underset{\textbf{v}_h \in X_h}{sup} \frac{\mid b(\textbf{v}_h,p-P_h)\mid}{\mid \textbf{v}_h \mid_{1,\Omega}}+\underset{\theta_h \in R_h}{sup}\frac{\frac{\rho_{hnf}}{\rho_{f}}\mid[a_3(T-T_h,\theta_h)+a_4(\textbf{u};T,\theta_h)-a_4(\textbf{u}_h;T_h,\theta_h)]\mid}{\mid \theta_h \mid_{1,\Omega}}
\end{equation}
Now we would like to estimate the non-linear terms 
\begin{equation}
\begin{split}
&\mid a_1(\textbf{u};\textbf{u},\textbf{v}_h)-a_1(\textbf{u}_h;\textbf{u}_h,\textbf{v}_h)\mid \\&= \mid a_(\textbf{u};\textbf{u}-\textbf{u}_h,\textbf{v}_h)+a_1(\textbf{u};\textbf{u}_h,\textbf{v}_h)-a_1(\textbf{u}_h;\textbf{u}_h,\textbf{v}_h)\mid\\
&\leq 2 N_1 R {\mid \textbf{u}-\textbf{u}_h \mid}_{1,\Omega} {\mid \textbf{v}_h \mid}_{1,\Omega} 
\end{split}
\end{equation}
\begin{equation}
\begin{split}
&\mid a_4(\textbf{u};T,\theta_h)-a_4(\textbf{u}_h;T_h,\theta_h) \mid \\&=\mid (\textbf{u} \cdot \nabla T)\theta_h-(\textbf{u}_h \cdot \nabla T_h)\theta_h \mid\\
&=\mid \big[ ( \textbf{u}-\textbf{u}_h) \cdot \nabla T+\textbf{u}_h \cdot \nabla (T-T_h)\big]\theta_h \mid\\
&\leq N_2\big[ ({\mid \textbf{u}-\textbf{u}_h \mid}_{1,\Omega} {\mid T \mid}_{1,\Omega} + {\mid \textbf{u}_h \mid}_{1,\Omega} {\mid T-T_h \mid}_{1,\Omega}) {\mid \theta_h \mid}_{1,\Omega}\big]\\
&\leq N_2R[{\mid \textbf{u}-\textbf{u}_h \mid}_{1,\Omega}+{\mid T-T_h\mid}_{1,\Omega}]{\mid \theta_h \mid}_{1,\Omega}.
\end{split}
 \end{equation}
 Estimation of the linear terms are as follows:\\
\begin{equation}
\mid a_0(\textbf{u}-\textbf{u}_h,\textbf{v}_h) \mid
\leq (\frac{\rho_f}{\rho_{hnf}})(\frac{\mu_{hnf}}{\mu_f}) Pr{ \mid \textbf{u}-\textbf{u}_h \mid}_{1,\Omega}{\mid \textbf{v}_h \mid}_{1,\Omega}
\end{equation}\vspace{2mm}\\
\begin{equation}
\mid a_3(T-T_h,\theta_h) \mid
\leq \frac{\alpha_{hnf}}{\alpha_f} {\mid T-T_h \mid}_{1,\Omega}{\mid \theta_h \mid}_{1,\Omega}
\end{equation}\vspace{2mm}\\
\begin{equation}
\mid a_2(j(T-T_h),\textbf{v}_h) \mid \leq \frac{(\rho \beta)_{hnf}}{\rho_{hnf}\beta_f}Pr Ra c^2 {\mid T-T_h \mid }_{1,\Omega}{\mid \textbf{v}_h \mid}_{1,\Omega}
\end{equation}
\vspace{3mm}\\
Therefore, using $(35)-(39)$ in equation $(34)$ we have \\
\begin{equation}
\|p_h-P_h\|_2 \leq \frac{\rho_{hnf}}{BL\rho_{f}}[((\frac{\rho_f}{\rho_{hnf}})(\frac{\mu_{hnf}}{\mu_f})Pr+2N_1R+N_2R)\|\nabla (\textbf{u}-\textbf{v}_h)\|_2
\nonumber \end{equation}
\begin{equation}
+(\frac{(\rho \beta)_{hnf}}{ \rho_{hnf}\beta_f}Pr Ra c^2  +N_2R+\frac{\alpha_{hnf}}{\alpha_f})\|\nabla (T-T_h)\|_2]+\frac{1}{BL}\|p-P_h\|_0.
 \end{equation}\vspace{3.5mm}\\
Using interpolation error estimate, (31) and (32) in inequality (40) we have\\
\begin{equation}
\begin{split}
&\|p_h-P_h\|_2 \\
&\leq \frac{\rho_{hnf}}{BL\rho_{f}}[((\frac{\rho_f}{\rho_{hnf}})(\frac{\mu_{hnf}}{\mu_f})Pr+2N_1R+N_2R) (c_1 h^2 \mid \textbf{u} \mid_{3}+c_2 h^2 \mid p \mid_{2}+c_3 h^2 \mid T \mid_{3})\\
&+(\frac{(\rho \beta)_{hnf}}{ \rho_{hnf}\beta_f}Pr Ra c^2 +N_2R+\frac{\alpha_{hnf}}{\alpha_f})(c_1 h^2 {\mid \textbf{u} \mid}_{3}+c_2 h^2 {\mid p \mid}_{2}+c_3 h^2 {\mid T \mid}_{3})]\\& +\frac{Ch^2}{BL}\mid p \mid _2.
\end{split}
 \end{equation}\vspace{3mm}\\
 Therefore,
 \begin{equation}
 \begin{split}
 &\|p-P_h\|_2\\
 &\leq \|p-P_h\|_2+\|p_h-P_h\|_2\\
 &\leq [((\frac{\rho_f}{\rho_{hnf}})(\frac{\mu_{hnf}}{\mu_f})Pr+2N_1R+N_2R) 
 +(\frac{(\rho \beta)_{hnf}}{ \rho_{hnf}\beta_f}Pr Ra c^2  +N_2R+\frac{\alpha_{hnf}}{\alpha_f})]\\
 &\times(c_1 h^2 {\mid \textbf{u} \mid}_{3}+c_2 h^2 {\mid p \mid}_{2}+c_3 h^2 {\mid T \mid}_{3})\times \frac{\rho_{hnf}}{BL\rho_{f}} +(1+{(BL)}^{-1})Ch^2\mid p \mid _2.
\end{split}
 \end{equation}\vspace{3.4mm}\\
 \texttt{Remark:} We have higher order of convergence with higher order of regularity assumption on the weak solutions\cite{SGS10}.
 .\vspace{7mm}\\
 \section{Numerical Experiments:}
To begin with both  the finite element scheme cum code validation tests  have been carried out based on the  benchmark problem from \cite{Gcp17}. Prior to carrying out the detailed simulations, the grid validation tests have been carried out both on H and L shaped enclosures. We employed an iterative strategy to develop our numerical scheme, primarily employing Cauchy linearization for nonlinear components.
The Newton quasi linearization approach is used to solve the nonlinear system produced by the finite element formulation.
While linear elements are used to estimate pressure using the $H^{1}(\Omega)$ norm, we have employed quadratic elements to approximate U, V, and T. We have calculated the numerical error to an accuracy up to a tolerance of $10^{-6}$. The terms "HNF1", "HNF2", and "HNF3" in the sections that follow refer to hybrid nanofluids having volume fractions of  $0.1\%$ $0.33\%$ and $1\%$ respectively.
\cite{CSM19} for information on the (Cu, Al2O3)/water hybrid nanofluid's thermophysical characteristics.

\subsection{Scheme validation:}
We have shown the comparison of Global Nussselt number$(Nu)$ calculated for the benchmark problem in \cite{Gcp17} consisting of the following geometry and boundary conditions:
square cavity with hot left vertical wall, cold right vertical wall and insulated horizontal top and bottom boundaries. Velocity is set to zero on all the boundaries. The obtained numerical result presented in Table 1 is in good agreement with the existing literature.
\begin{table}[!h]
	\centering
	\scalebox{1.5}{}
	
	\begin{tabular}{|l|c|c|c|}
		\hline
		&\multicolumn{3}{c}{|Nusselt No.|} \\
		\cline{2-4}
		\newline
		Paper & $Ra=10$ & $Ra=100$ & $Ra=10^3$ \\
		\hline
		\hline
		Walker et al. \cite{Who78}  & -  & 3.10 & 12.96  \\
		\hline
		Bejan \cite{Bej79}    & -  & 4.20 & 15.80 \\
		\hline
		Beckerman et al. \cite{Bvr86}   & - & 3.11 & - \\ 
		\hline
		Moya et al. \cite{Mrs87} & 1.065 & 2.80 & - \\
		\hline
		Manole and Lage \cite{Mal92} & -& 3.12 & 13.64 \\
		\hline
		Baytas and Pop \cite{Bap02} & 1.079 & 3.16 & 14.06 \\
		\hline
		Giovani et al. \cite{Gcp17}   & 1.078 & 3.12 & 15.60 \\
		\hline
		Our work   & 1.0791 & 3.13181 & 14.8515 \\
		\hline
	\end{tabular}
	\caption{Comparison of Nusselt no. value with standard literature for validation of the numerical scheme:}
\end{table}

\subsection{Grid validation test:}
A grid validation test, for $Pr=50,Ra=10^7,\phi=1\%$, is done based on Nusselt number$(Nu)$, $\psi_{max}$, $\psi_{min}$ is presented in Table 2. This experiment is done on an L shape domain with the following  set up: left vertical wall is kept hot and right vertical is cold and rest of the walls are kept adiabatic. We have picked up the grid size $100 \times 100$ as we observe from Table 2 that starting from grid size $100 \times 100$ the change in  Nusselt no., $\psi_{max}$, $\psi_{min}$ is negligible. So, further numerical computations on L shaped domain are performed with the grid size $100 \times 100$.
\begin{table}[h]
	\centering
	\scalebox{5}{}
	\begin{tabular}{l|c|c|c|r}
		\hline
		
		Grid size & Global Nu & $\psi_{max}$ &$\psi_{min}$ \\
		\hline
		\hline
	
		$40 \times 40$  & 25.1091737   & 232.351 &-7.19562  \\
		\hline
		$ 60 \times 60$  & 24.7312602   & 233.375 & -6.76847 \\
		\hline
		$80 \times 80$  & 24.6504471  & 234.282 & -7.29636 \\
		\hline 
		$90 \times 90$ & 24.6546271 & 234.39 & -7.23942\\
		\hline
		$100 \times 100$ & 24.6223662 & 234.6 &-7.22445\\
		\hline
		$110 \times 110$ & 24.604896 & 234.68 &-7.32931\\
		\hline
		$120 \times 120$ & 24.6102549 & 234.767 &-7.25268\\
		\hline
		
	\end{tabular}
	\caption{Grid test on L shape domain}
\end{table}
Similarly, we have carried out the grid validation test on H shaped cavity with the following parametric set up like $Ra=10^5, Pr=10$ and $\phi=1\%$. One can observe from table 3 that starting from the grid size $60 \times 60$ the variation in obtained data like global Nusselt number$(Nu)$, $\psi_{max}$, $\psi_{min}$ is negligible. So, further computations on H shaped cavity are carried out on a $64 \times 64$ grid.
\begin{table}[h]
	\centering
	\scalebox{5}{}
	\begin{tabular}{l|c|c|c|r}
		\hline
		
		Grid size & Global Nu & $\psi_{max}$ &$\psi_{min}$ \\
		\hline
		\hline
		
		$20 \times 20$ & 7.1901 & 13.0793 & -0.0977427\\
		\hline
		$30 \times 30$ & 7.14068 & 13.0766 & -0.117362 \\
		\hline 
		$40 \times 40$  & 7.12643 & 13.0742 & -0.117207 \\
		\hline
		$50 \times 50$  & 7.12087 & 13.0754 & -0.121484 \\
		\hline
		$60 \times 60$  & 7.11831 &  13.0756&  -0.119795\\
		\hline
		$64 \times 64$  & 7.11780 & 13.0752  &  -0.119638\\
	   \hline
		
	\end{tabular}
	\caption{Grid test on H shape domain}
\end{table}
 
  \section{H Shape Domain:}
  Here we have considered the following H shaped domain with $64X64$ grid size shown in Fig.1. The numerical experiment is done with the following boundary conditions:\\
  Extreme left vertical wall is hot and is maintained at $T=1.$\\
  Extreme right vertical wall is cold and is maintained $T=0.$\\
  Rest of the walls of the domain are adiabatic.\\
  $U=V=0 $ throughout the boundary.
 \begin{figure}[h]
  \centering
  \includegraphics[width=85mm]{ 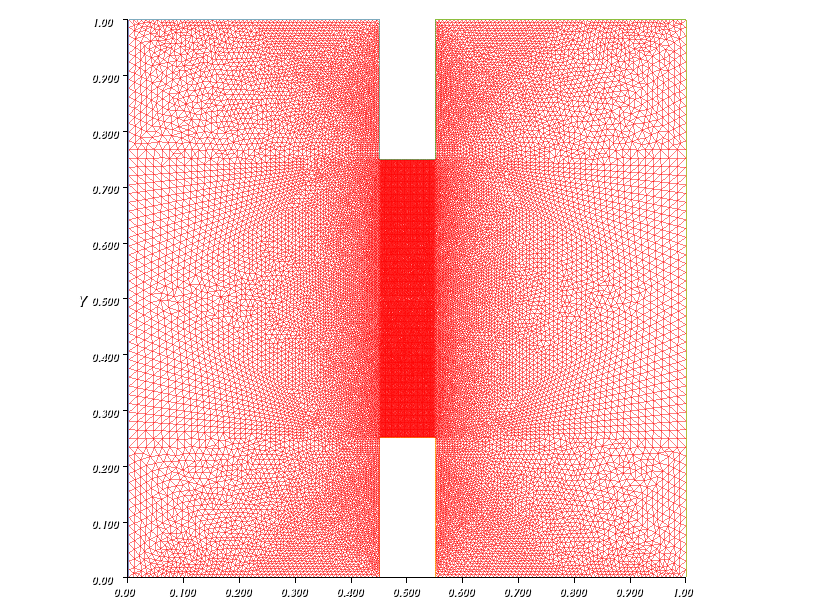}
    \caption{ H Shape Domain}
    \label{fig:height .25 unit,breadth 1 unit}
  \end{figure}
  \subsection{ Comparison of our result and with existing literature for code validation:}

\begin{table}[!h]
	\centering
	\scalebox{1.5}{}
	
	\begin{tabular}{|l|c|c|c|}
		\hline
		&\multicolumn{3}{c}{|Nusselt No.|} \\
		\cline{2-3}
		\newline
		Paper & $Ra=10^3$  & $Ra=10^4$ \\
		\hline
		\hline
		CFX,2001 \cite{CFX01}   & 0.751  & 1.267  \\
		\hline
		Skerget et al. \cite{LSH03} & 0.753  & 1.231 \\
		\hline
		Present Study    & 0.75371 & 1.33099   \\
		\hline
	\end{tabular}
	\caption{Comparison of Nusselt no. value with standard literature for validation of the numerical scheme:}
\end{table}
Because they closely mirror those in the body of existing literature, the numerical results in Table 4 justify our numerical scheme on the H shape domain.
Additionally, we have shown a few plots that resemble those in the literature cited in \cite{LSH03}.
To match those in \cite{LSH03}  who also utilised the HSV colour scheme, we have adopted the HSV colour scheme in the case of figures 2-7.
However, we adopted an RGB colour scheme throughout the entire manuscript to display the results of our own experiments using streamlines and isotherms. 
\begin{figure}[!h]
  \centering
  \begin{minipage}[b]{0.47\textwidth}
    \includegraphics[width=\textwidth]{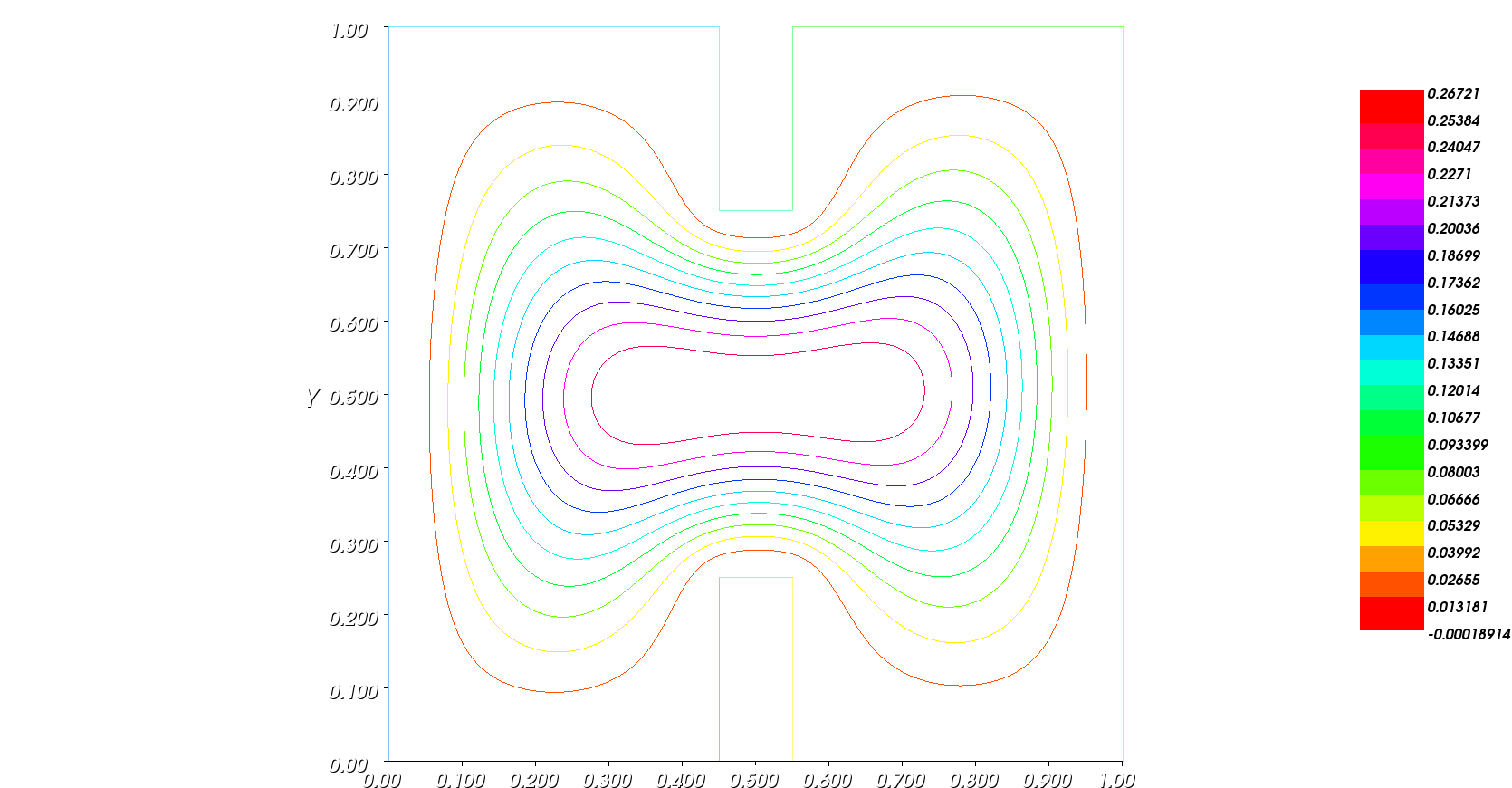}
    \caption{Streamlines for CF $Ra=10^3$ }
  \end{minipage}
  \hfill
  \begin{minipage}[b]{0.47\textwidth}
    \includegraphics[width=\textwidth]{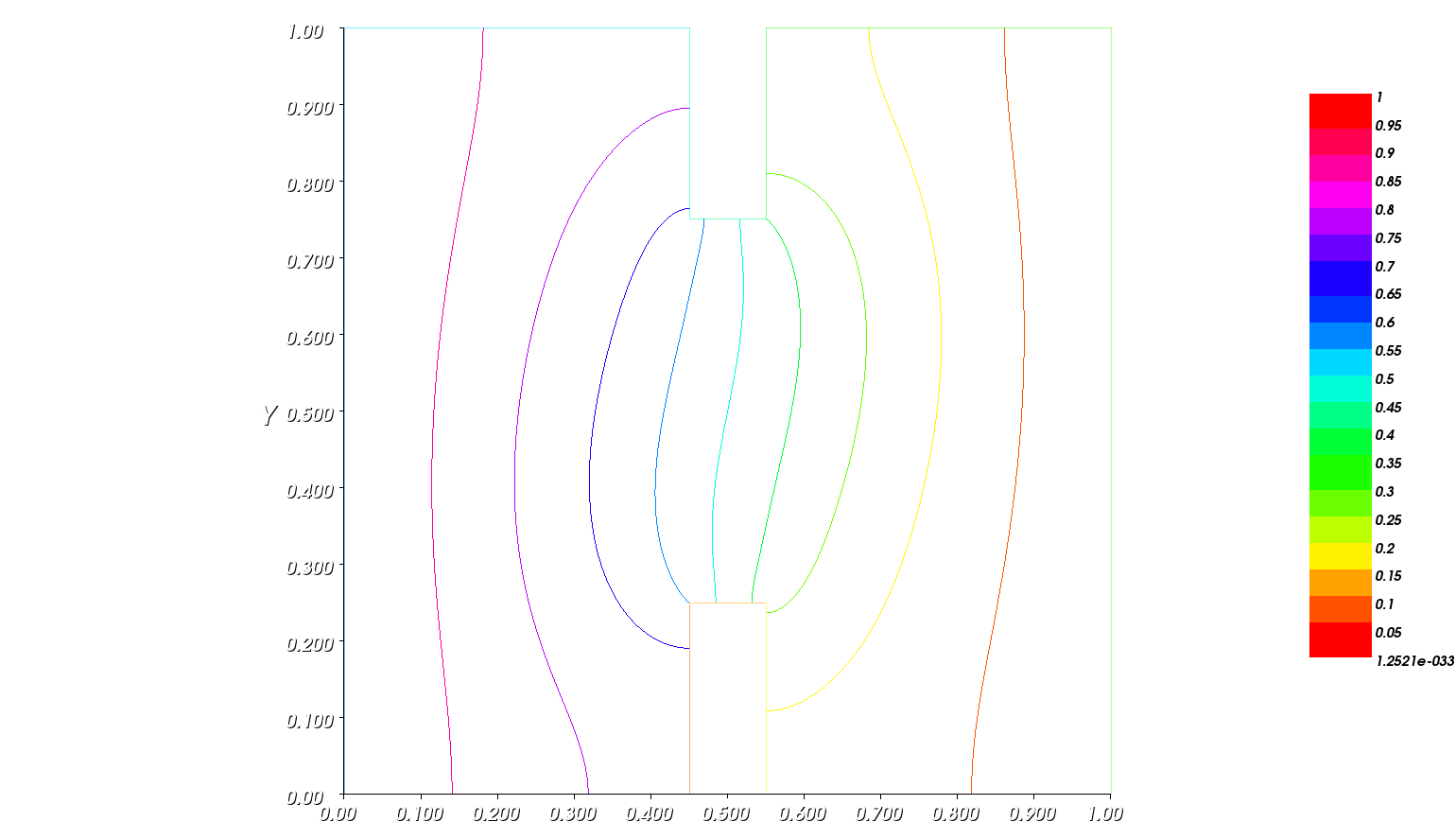}
    \caption{Isotherms for CF $Ra=10^3$}
  \end{minipage}
  \end{figure}
  \begin{figure}[!h]
  \centering
  \begin{minipage}[b]{0.47\textwidth}
    \includegraphics[width=\textwidth]{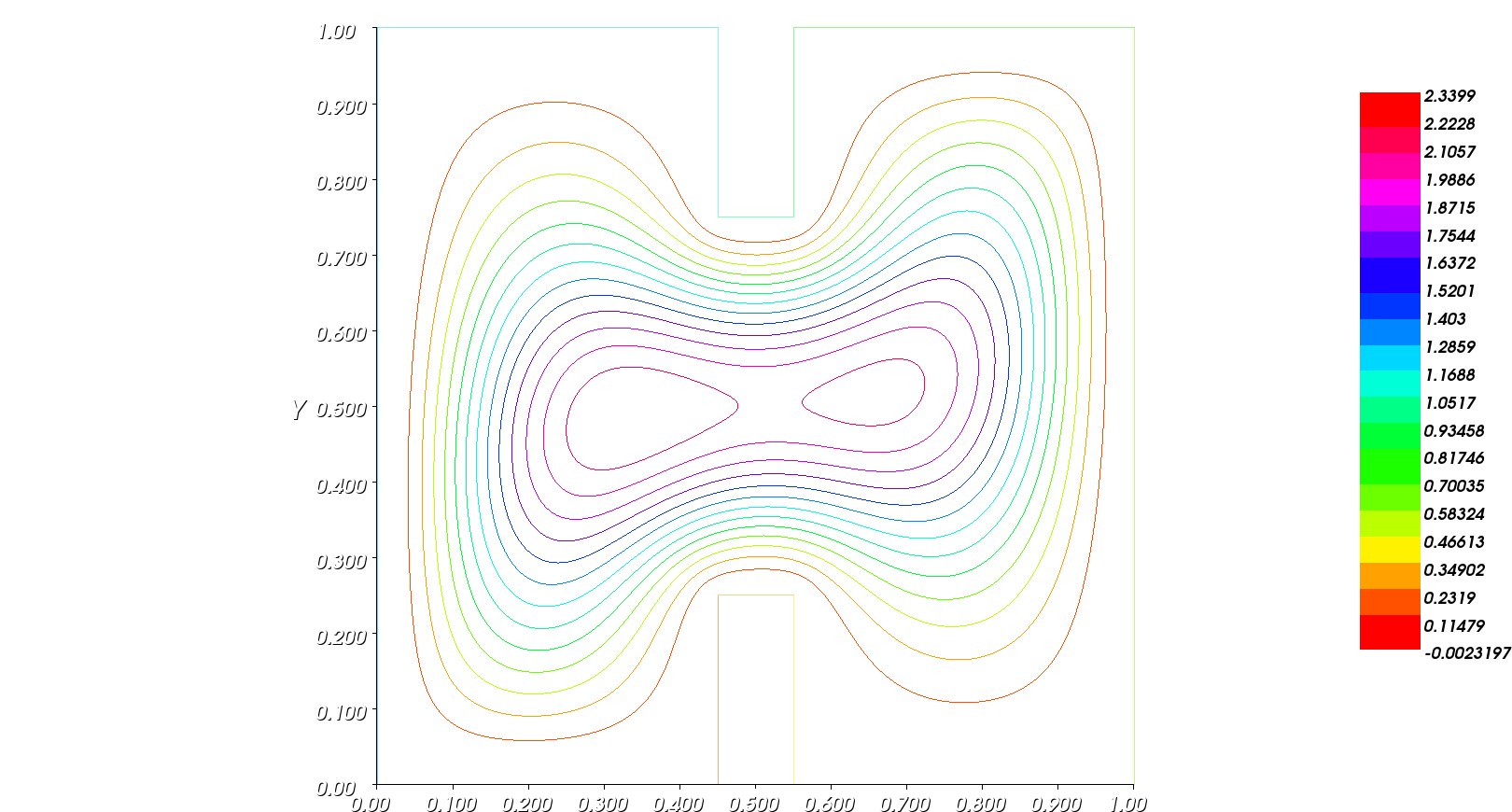}
    \caption{Streamlines for CF $Ra=10^4$ }
  \end{minipage}
  \hfill
  \begin{minipage}[b]{0.47\textwidth}
    \includegraphics[width=\textwidth]{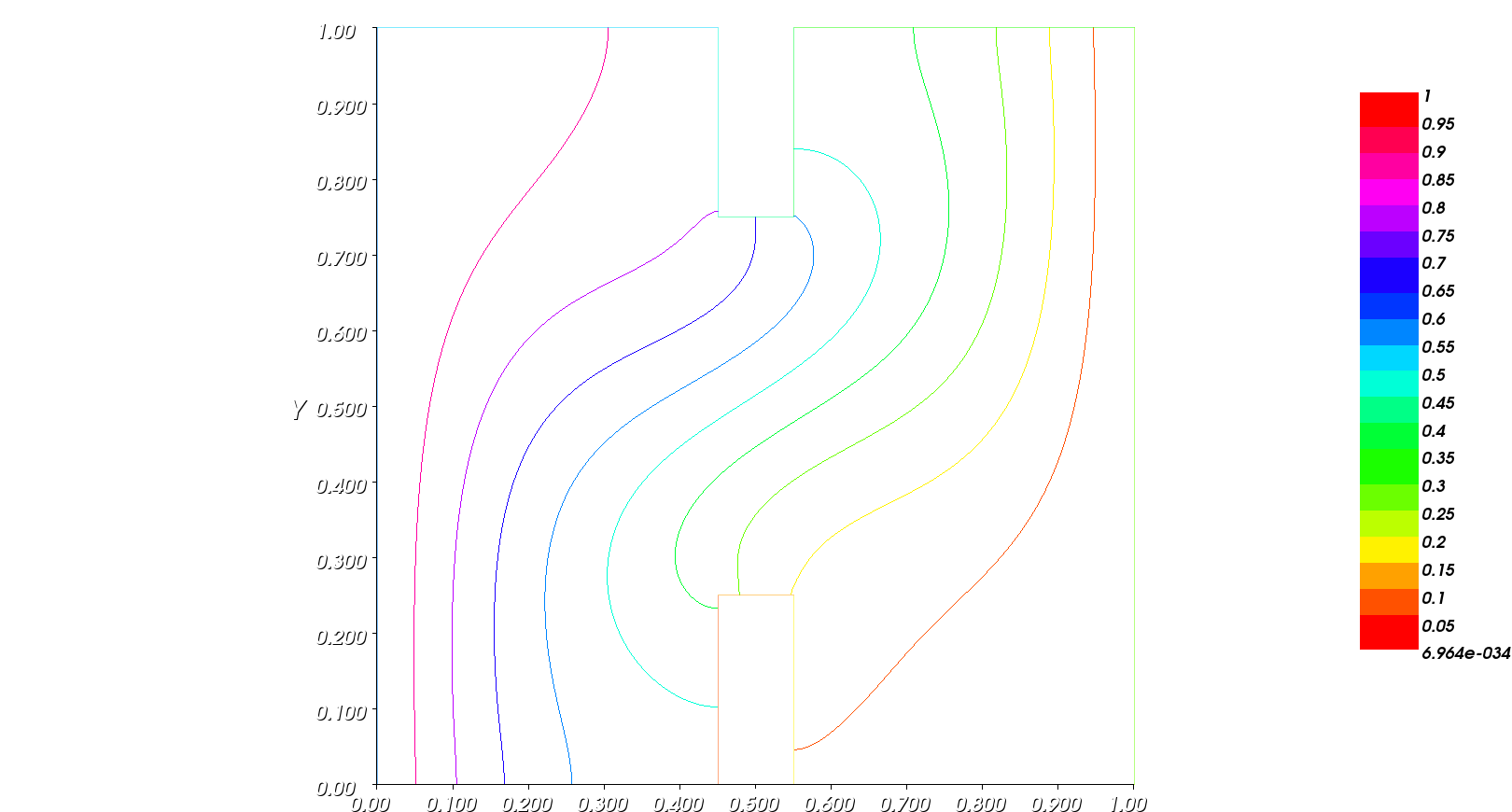}
    \caption{Isotherms for CF $Ra=10^4$}
  \end{minipage}
  \end{figure}
  \begin{figure}[!t]
  \centering
  \begin{minipage}[b]{0.47\textwidth}
    \includegraphics[width=\textwidth]{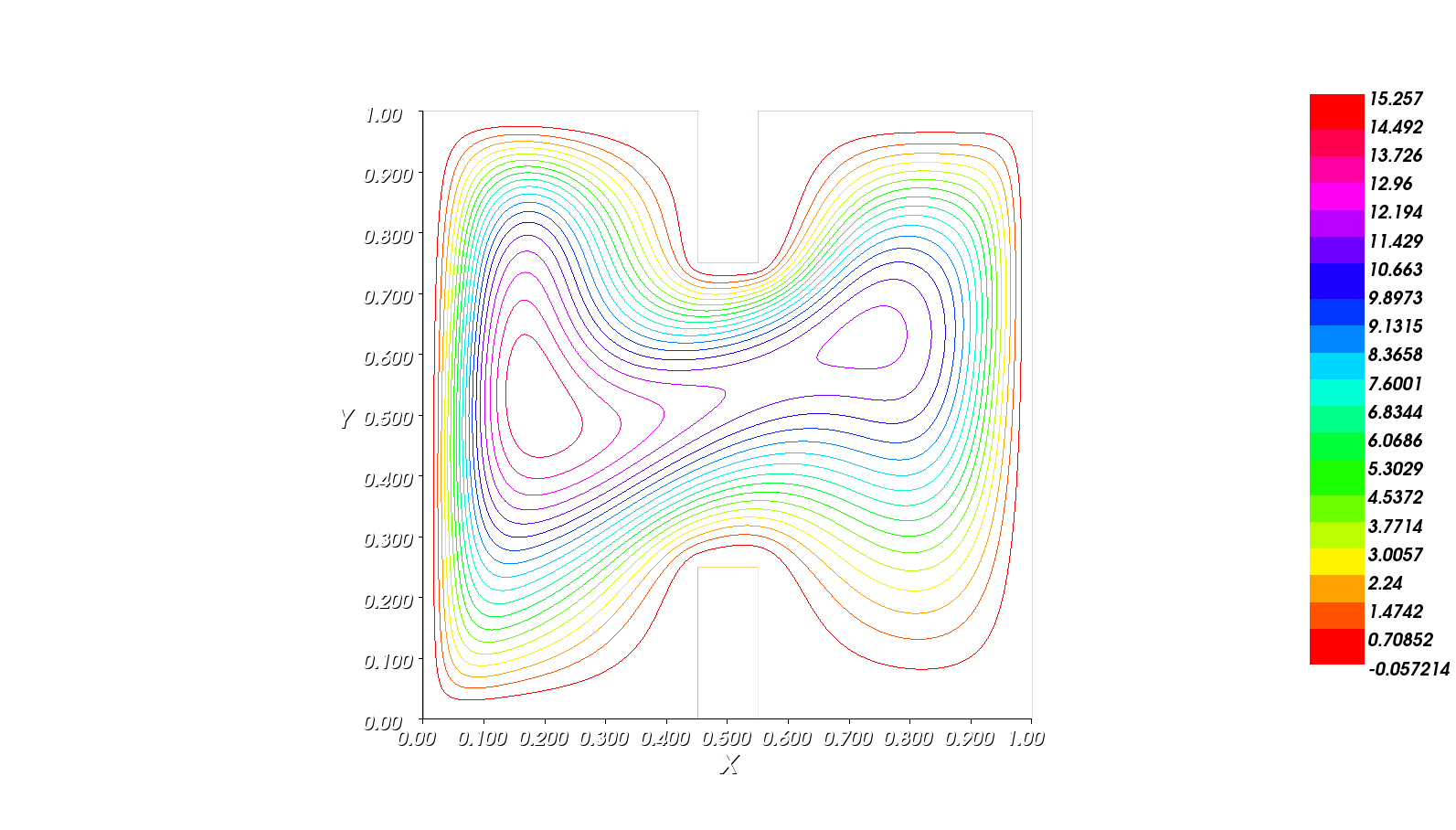}
    \caption{Streamlines for CF for $Ra=3.5X10^5$ }
  \end{minipage}
  \hfill
  \begin{minipage}[b]{0.47\textwidth}
    \includegraphics[width=\textwidth]{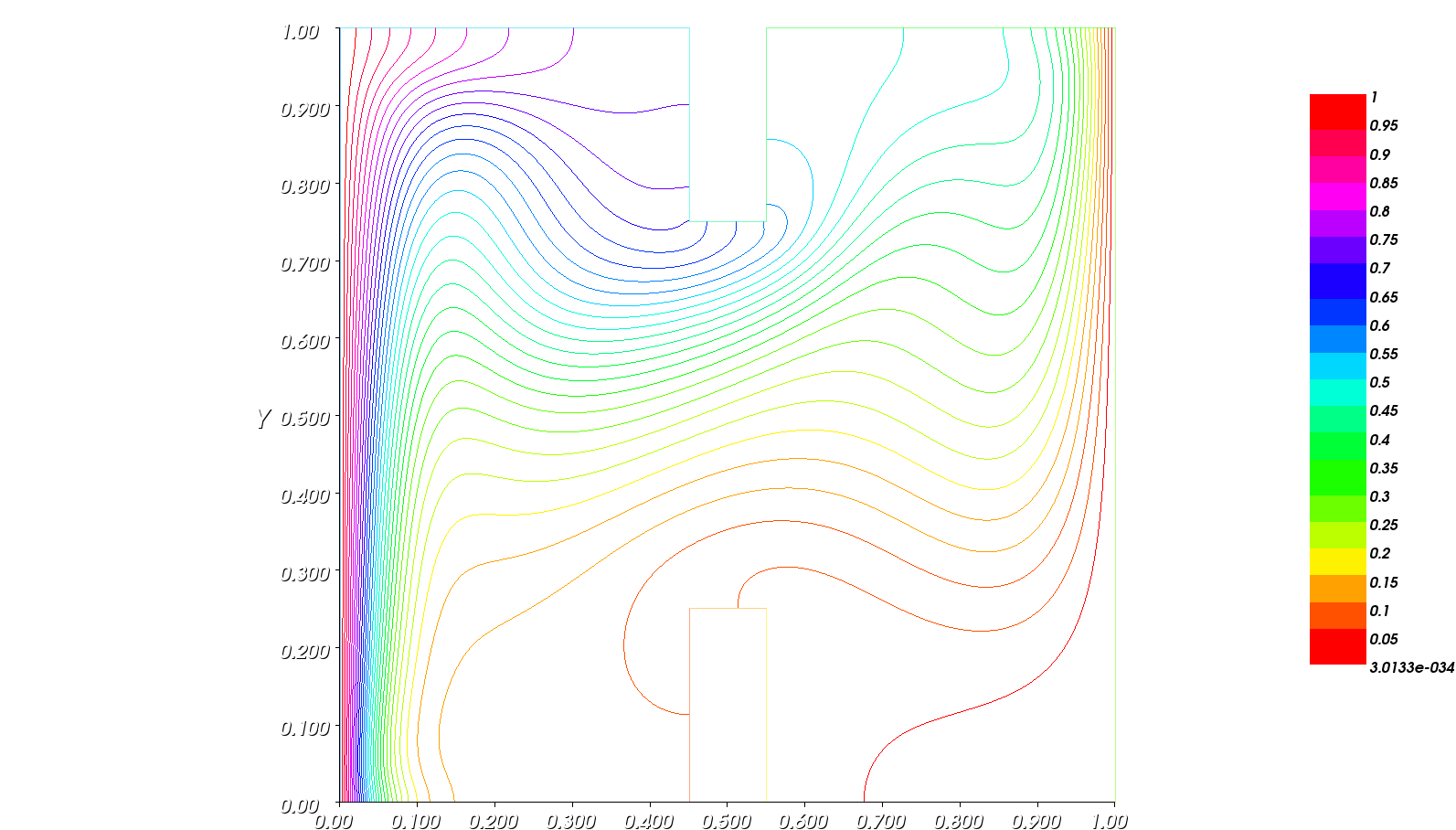}
    \caption{Isotherms for CF $Ra=3.5X10^5$}
  \end{minipage}
  \end{figure}\vspace{3mm}\\
  \clearpage
  \newpage
\subsection{The effect of input parameters $Ra,~\& ~\phi$, for fixed value of $Pr$, on Nu calculated on H shape enclosure:} 
 
\begin{table}[h]
	\centering
	\scalebox{5}{}
	\begin{tabular}{l|c|c|c|r}
		\hline
		
		$Pr=1,Ra$ &  $CF ~Nu$ & $HNF ~Nu$ & $HNF1 ~Nu$ & $HNF2~ Nu$\\
		\hline
		\hline
		 $Ra=1$  & 0.744341  & 0.752695104 & 0.766173385 & 0.797778126\\
		\hline 
		$ Ra=10^2$  & 0.744436 & 0.752788136 & 0.766264995 & 0.797864941\\
		\hline 
		$Ra=10^3$  & 0.75371  & 0.761942743 & 0.775216059 & 0.806344953\\
		\hline 
		$Ra=10^4$ & 1.33099 & 1.33590713 & 1.34285514 & 1.35921343 \\
		\hline
		$Ra=10^5$ & 4.9713 & 4.99457439 & 5.02917823 & 5.1023477\\
	\hline
		
	\end{tabular}
	\caption{Comparison of Nusselt no. for different volume fraction on H shape domain}
\end{table}
 We have presented a few selected plots below to explain the impact of the parameter $Ra$ on the fluid flow and heat transfer for the fixed value of  $\phi=1\%,Pr=1.$
 \begin{figure}[!h]
  \centering
  \begin{minipage}[b]{0.47\textwidth}
    \includegraphics[width=\textwidth]{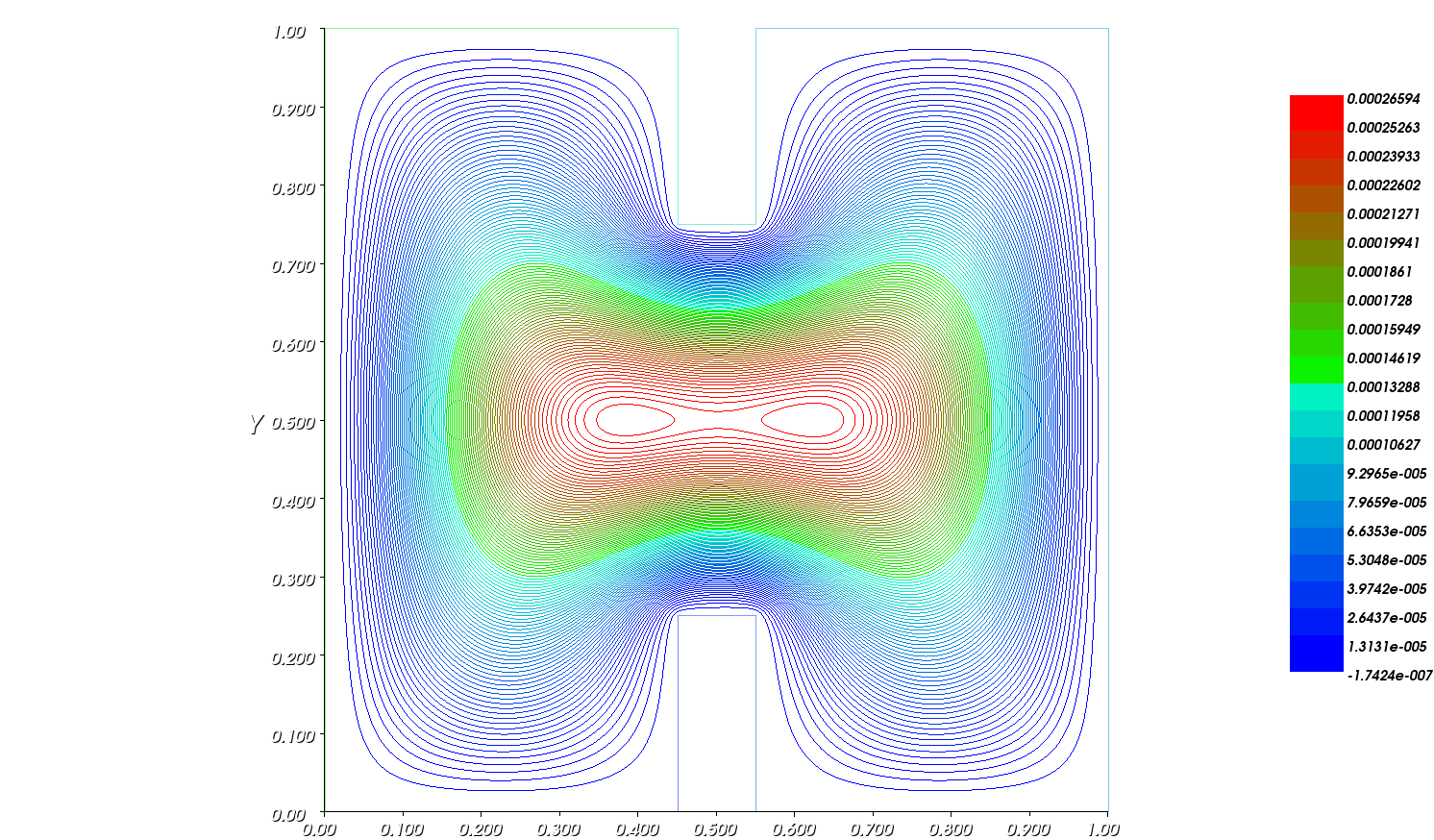}
    \caption{streamlines $Ra=1$  }
  \end{minipage}
  \hfill
  \begin{minipage}[b]{0.47\textwidth}
    \includegraphics[width=\textwidth]{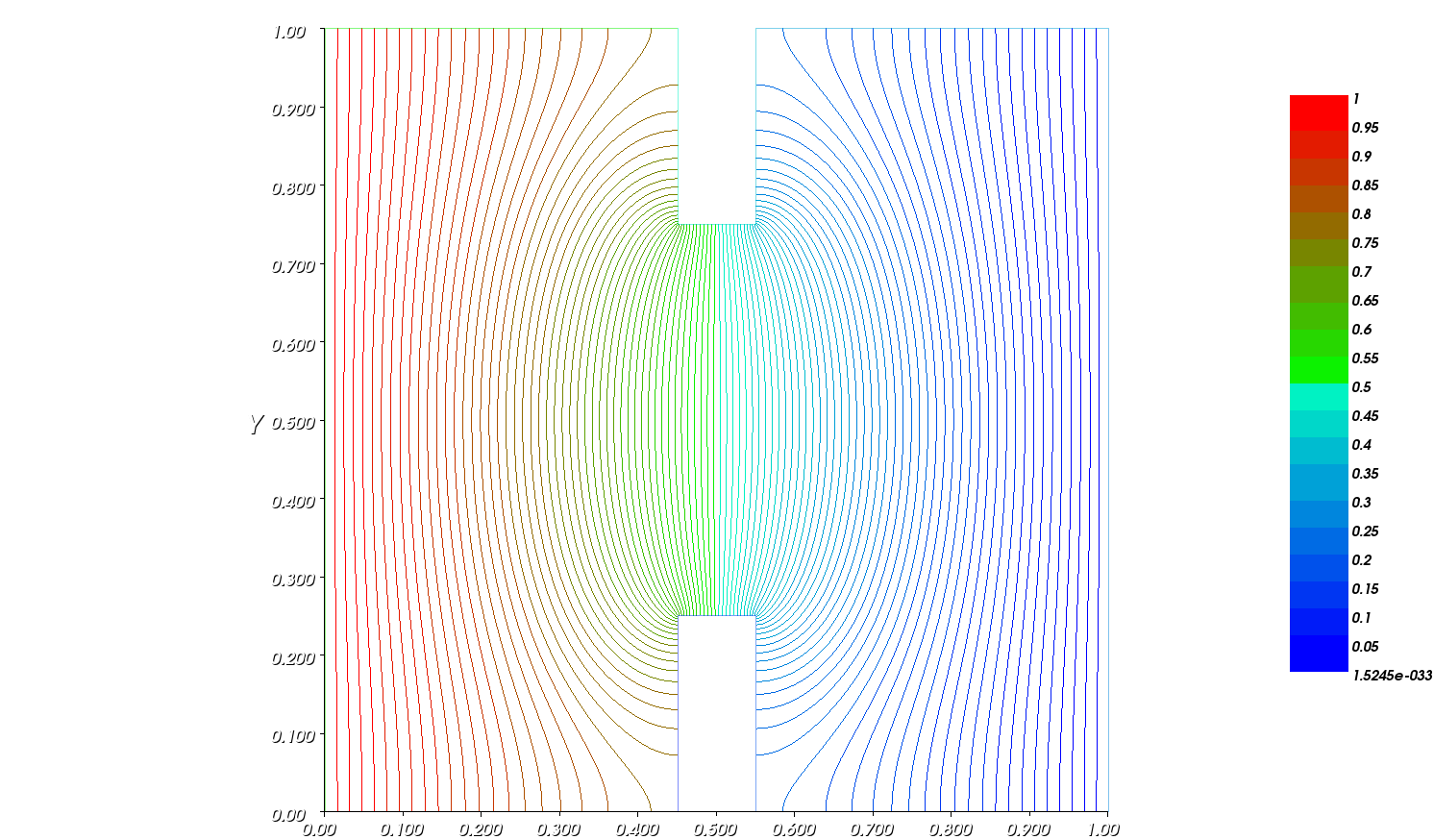}
    \caption{Isotherms $Ra=1$}
  \end{minipage}
  \end{figure}
 
 \begin{figure}[!h]
  \centering
  \begin{minipage}[b]{0.45\textwidth}
    \includegraphics[width=\textwidth]{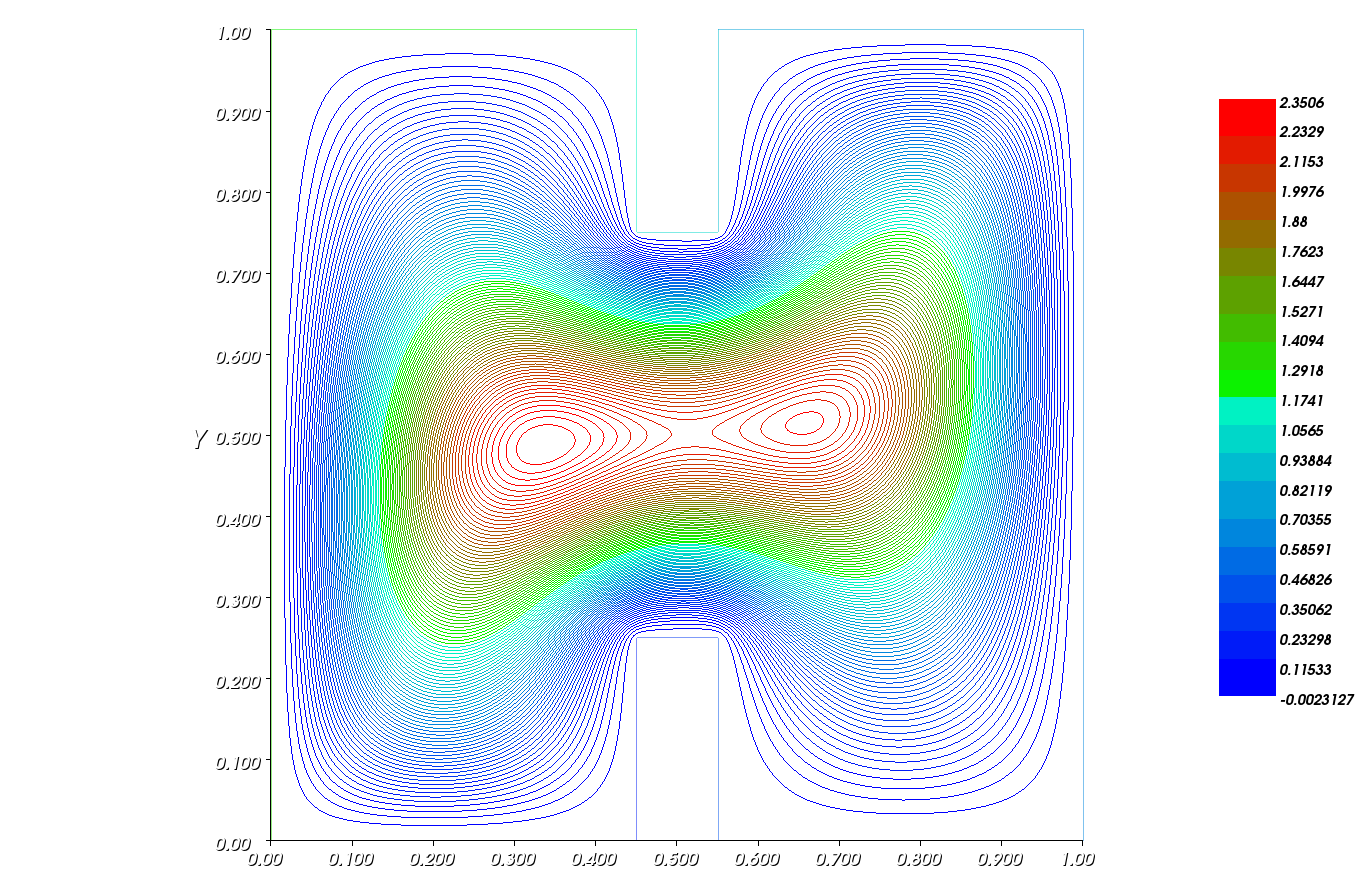}
    \caption{streamlines $Ra=10^4$  }
  \end{minipage}
  \hfill
  \begin{minipage}[b]{0.45\textwidth}
    \includegraphics[width=\textwidth]{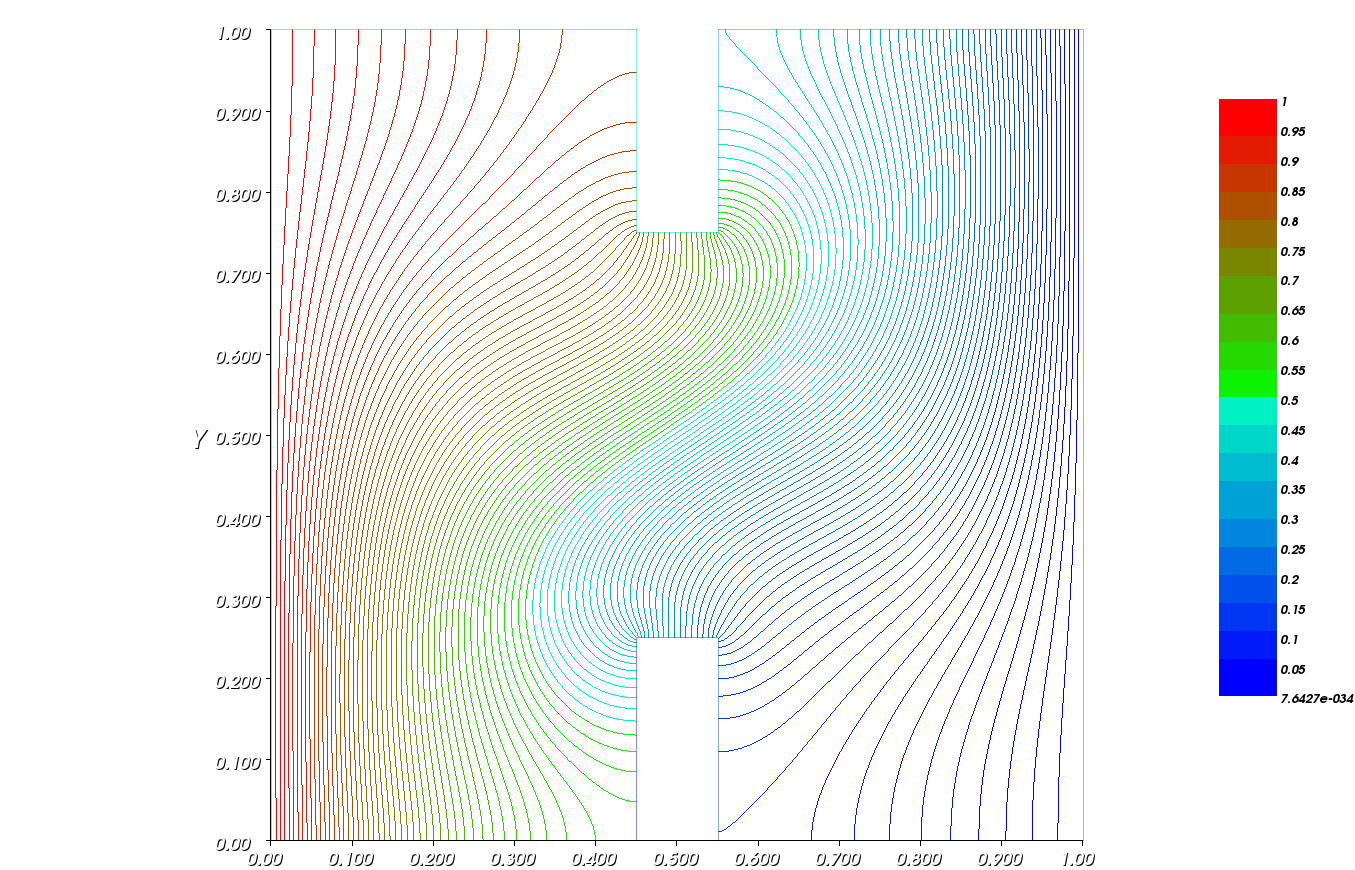}
    \caption{Isotherms $Ra=10^4$ }
  \end{minipage}
  \end{figure}
\clearpage
\newpage
   \begin{figure}[!h]
  \centering
  \begin{minipage}[b]{0.45\textwidth}
    \includegraphics[width=\textwidth]{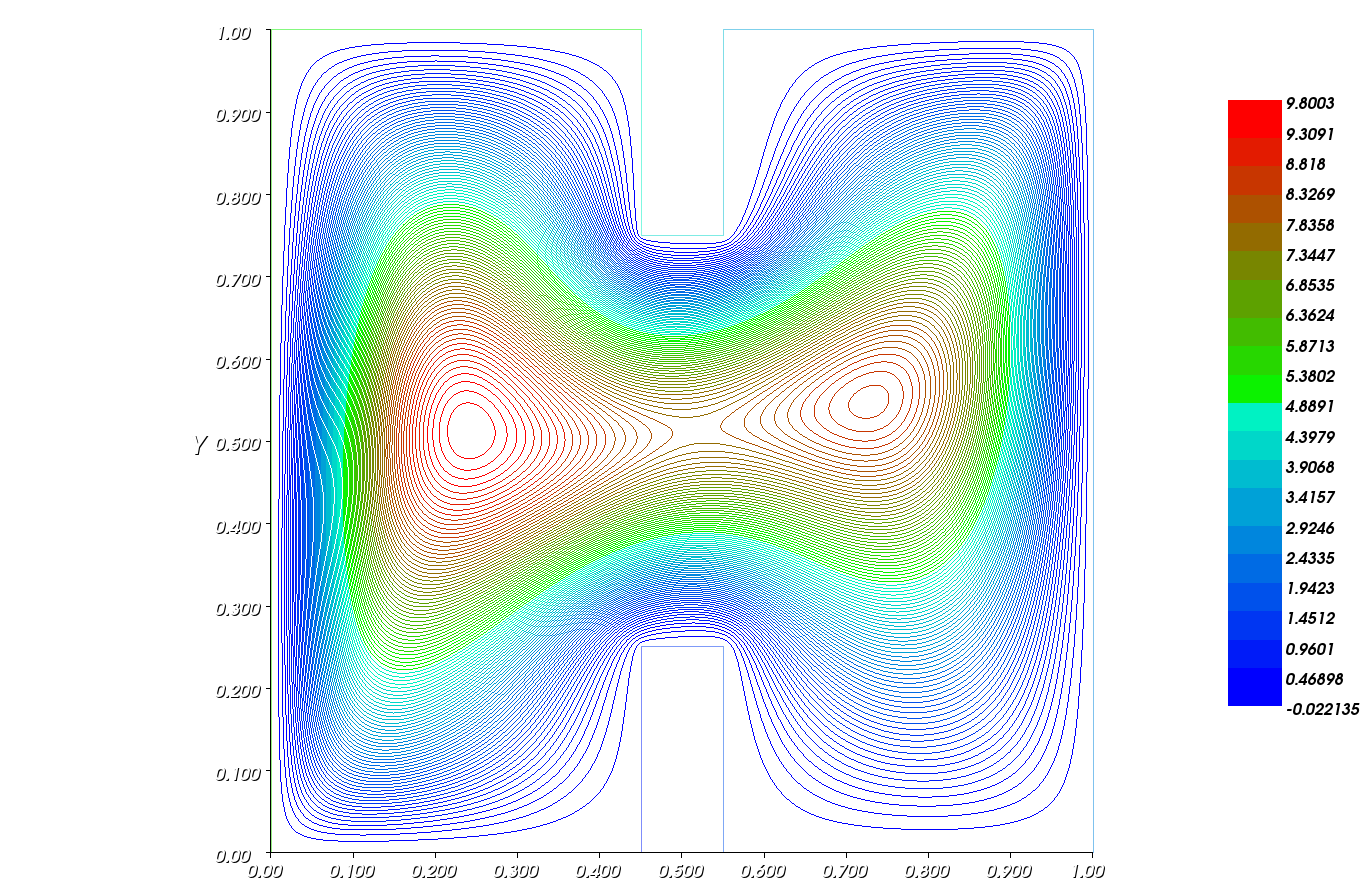}
    \caption{streamlines $Ra=10^5$ }
  \end{minipage}
  \hfill
  \begin{minipage}[b]{0.45\textwidth}
    \includegraphics[width=\textwidth]{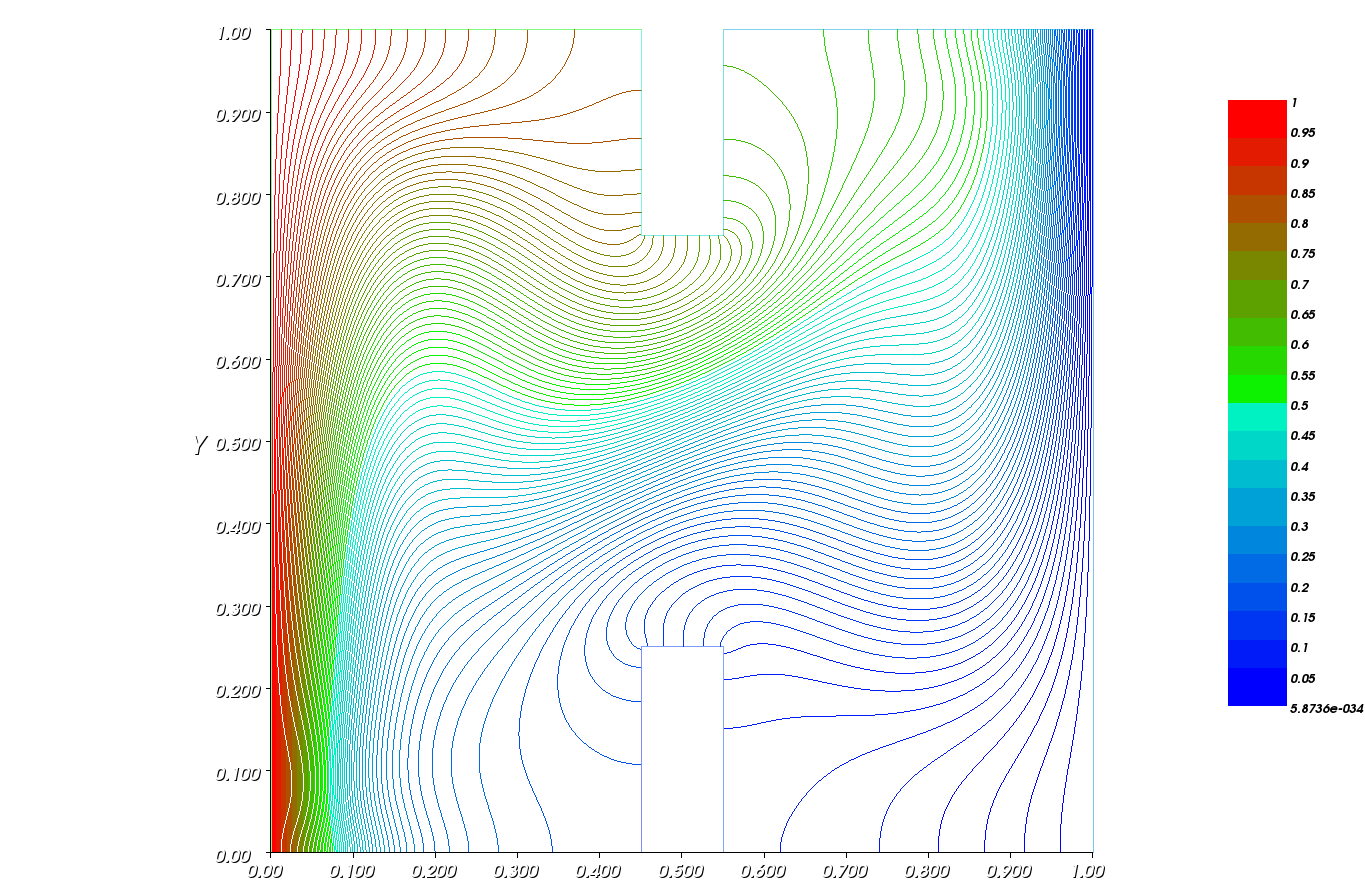}
    \caption{Isotherms $Ra=10^5$}
  \end{minipage}
  
  \end{figure}
   
\clearpage
\newpage
\subsection{Effect of input parameters $Pr,~\& \phi$, for fixed value of $Ra$, on Nu:}
  \clearpage
 \newpage
\begin{table}[h]
	\centering
	\scalebox{5}{}
	\begin{tabular}{l|c|c|c|r}
		\hline
		
		Pr, $(Ra=10^4)$ & CF Nu & HNF1 Nu & HNF2 Nu & HNF3 Nu \\
		\hline
		$Pr=1$  & 1.33099 &    1.33590713 & 1.3432257  & 1.35921343\\
		\hline
		$ Pr=5$  & 1.56383 &    1.5668908 & 1.57096522  & 1.57900565  \\
		\hline 
		$Pr=10$  &1.73341 &    1.73458199 & 1.73539058 & 1.73542286  \\
			\hline
		
	\end{tabular}
	\caption{Comparison of Nusselt no. for different volume fraction and Pr on H shape domain}
\end{table}
  

 We have presented a few plots to show the effect of the parameter $Pr$ on the flow pattern and temperature distributions for $\phi=0.33\%$.
\begin{figure}[!h]
  \centering
  \begin{minipage}[b]{0.32\textwidth}
    \includegraphics[width=\textwidth]{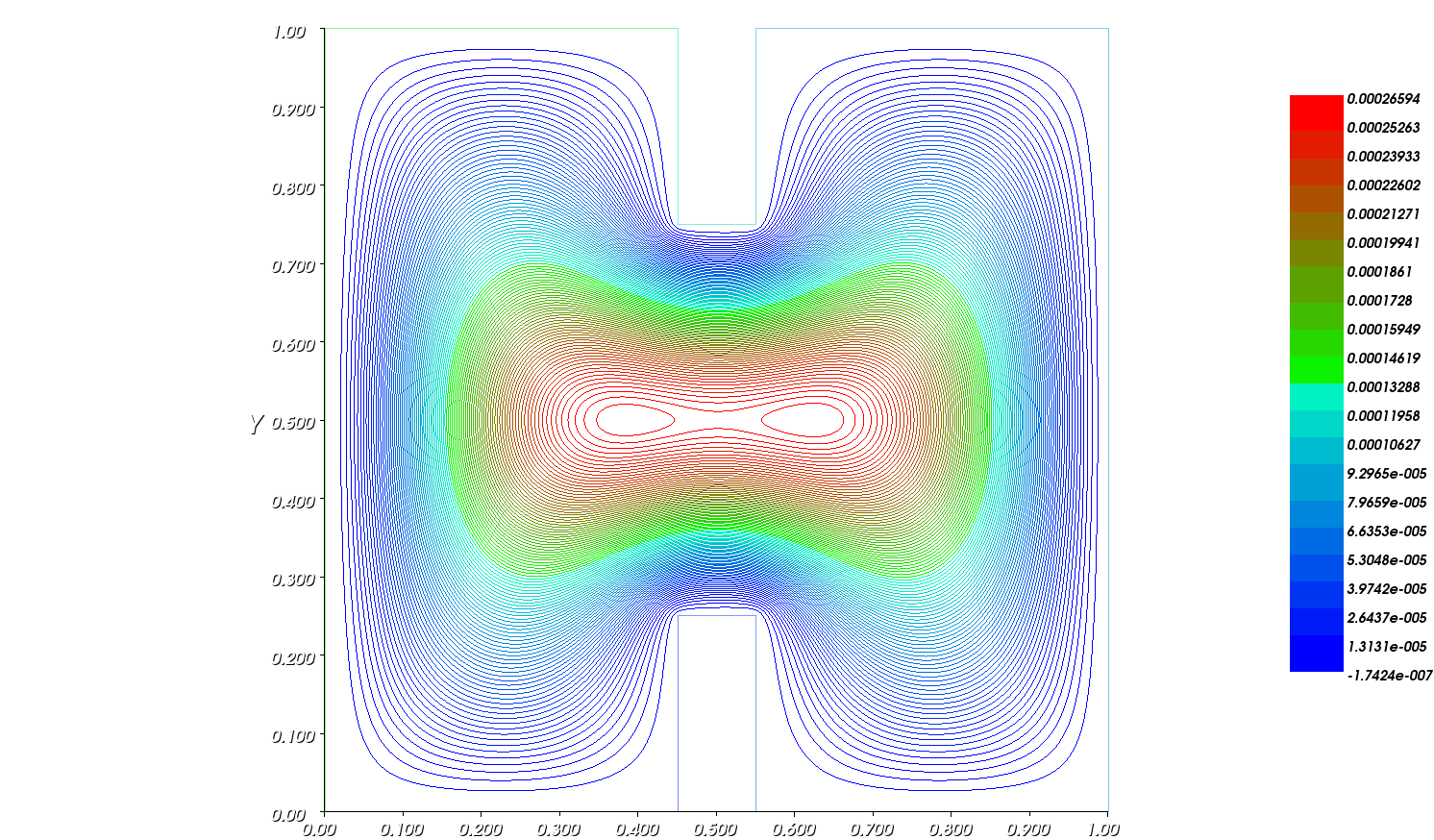}
    \caption{streamlines $Pr=1$ }
  \end{minipage}
  \hfill
  \begin{minipage}[b]{0.32\textwidth}
    \includegraphics[width=\textwidth]{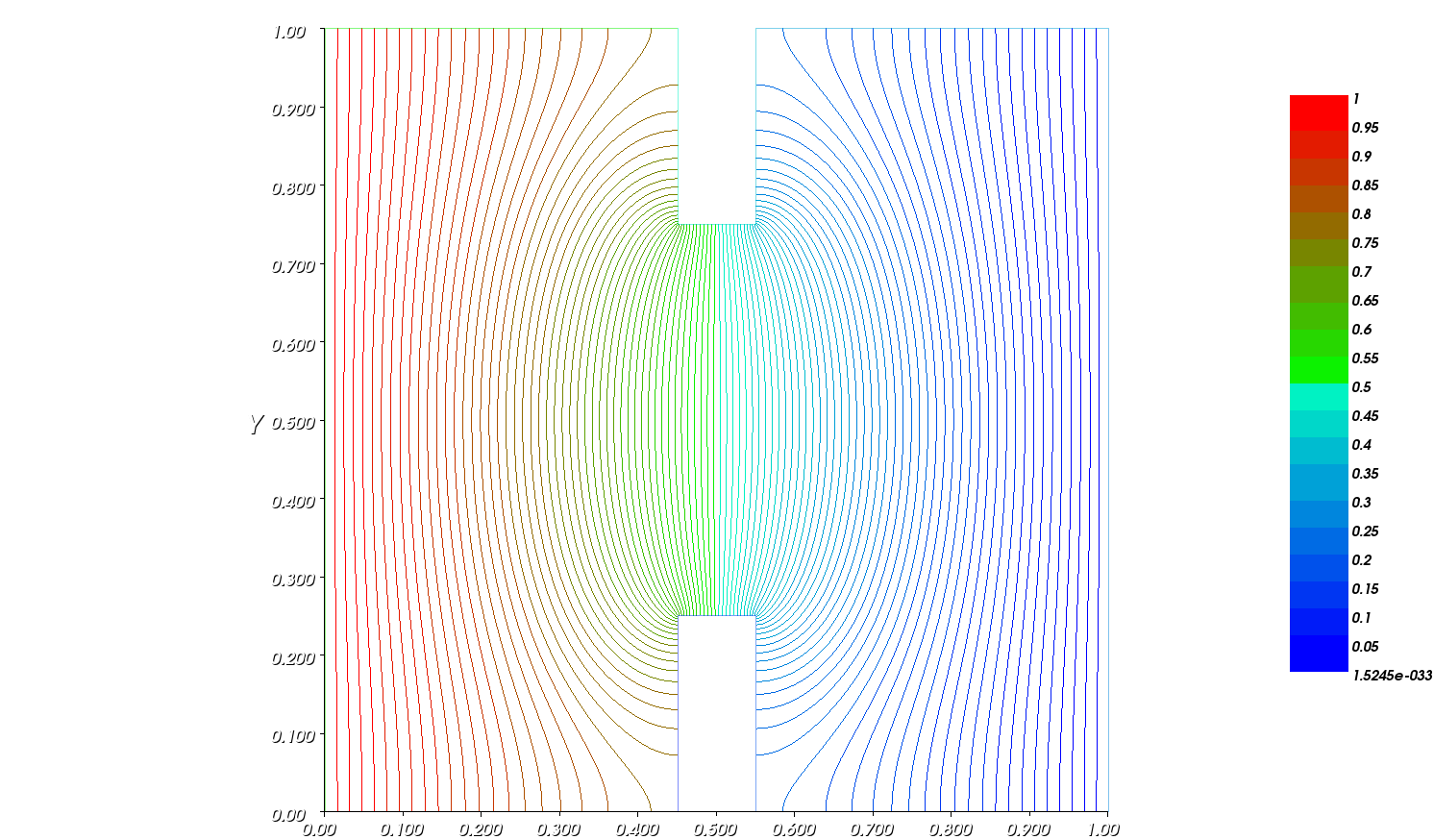}
    \caption{Isotherms $Pr=1$}
  \end{minipage}
  \end{figure}
   \begin{figure}[!h]
  \centering
  \begin{minipage}[b]{0.32\textwidth}
    \includegraphics[width=\textwidth]{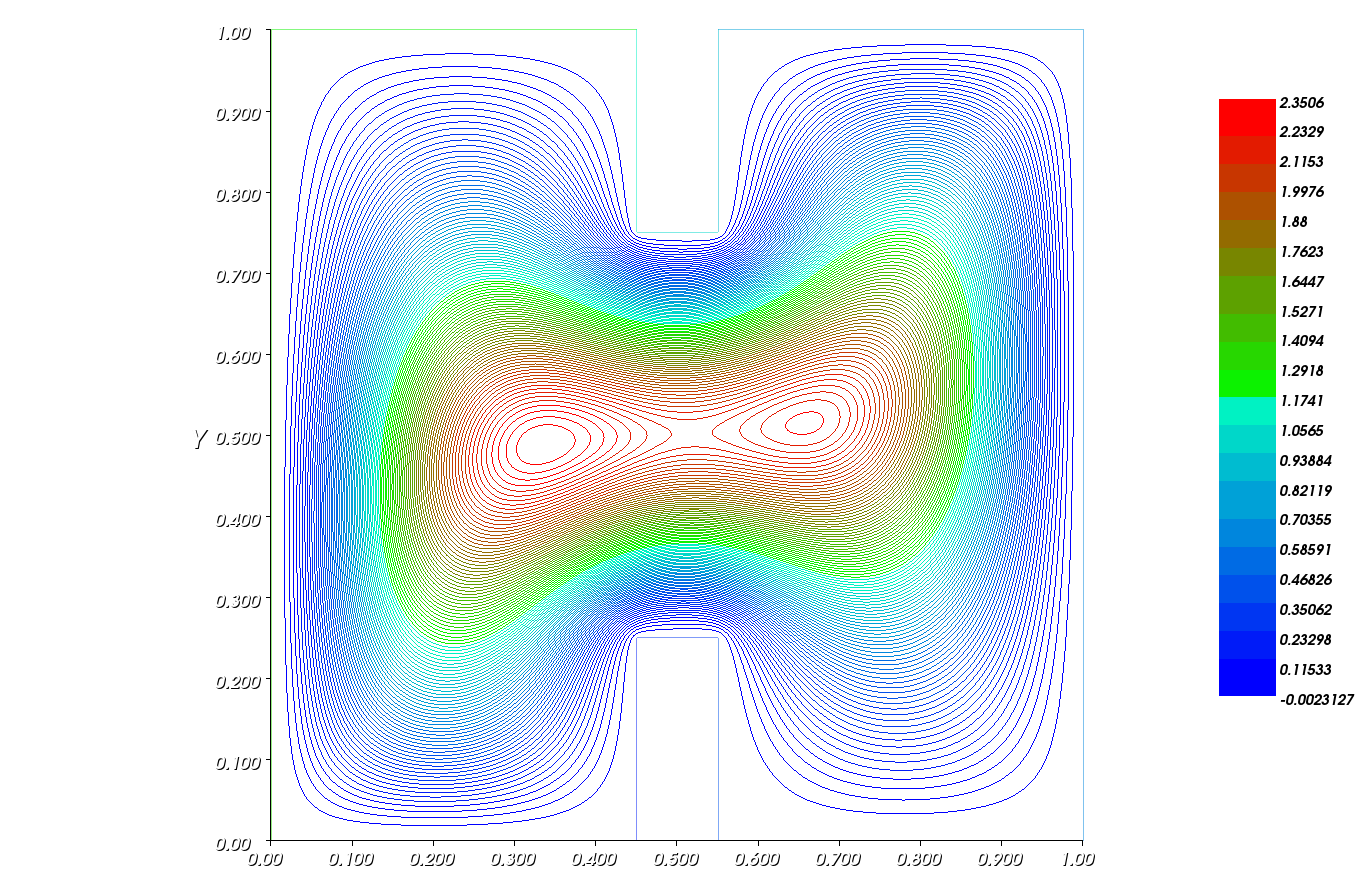}
    \caption{streamlines $Pr=10$ }
  \end{minipage}
  \hfill
  \begin{minipage}[b]{0.32\textwidth}
    \includegraphics[width=\textwidth]{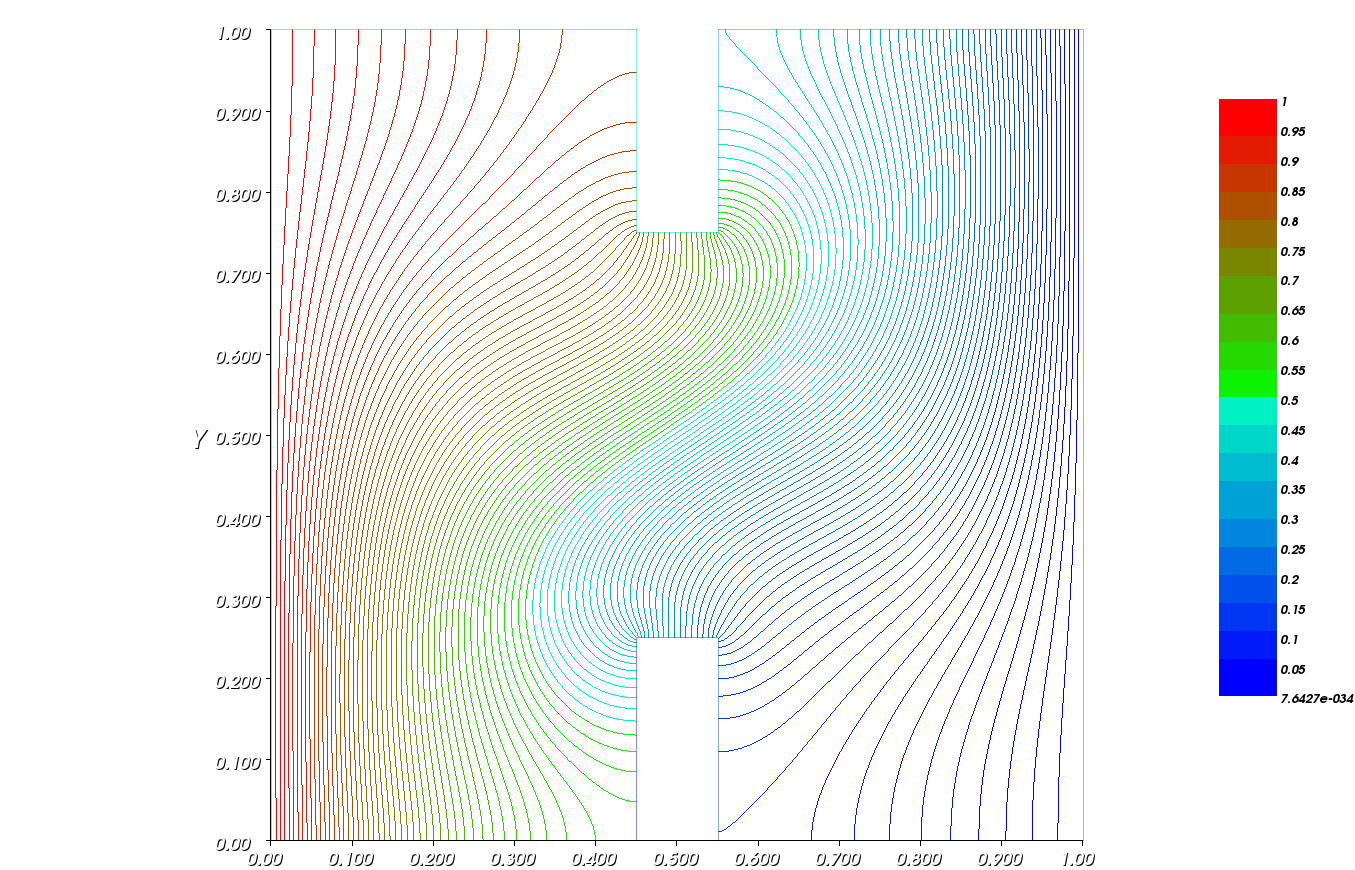}
    \caption{Isotherms $Pr=10$}
  \end{minipage}
  
  \end{figure}

\subsection{Results and discussions on H shape domain:}
One can observe, from fig. 8, fig.10, fig.12 depicting the streamlines, the flow pattern forms a butterfly like shape inside the domain with a dumbbell like shape situated in the center. When Ra is low, like $Ra=1$, the main circulation eye is located in the middle of the dumbbell, however when Ra is high, like $Ra=10^5$, two distinct, clearly visible circulation eyes are evident.
Additionally, as Ra is increased from $1$ to $10^5$, the circulation's eyes began to wander towards the direction of vertical walls spaced apart from one another.
Additionally, we can see the development of multicellular patterns that are co-rotating and symmetric in nature.
Also, when the value of $Ra$ increases from $Ra=1$ to $Ra=10^4$, more of these patterns are generated and are travelling in the direction of the domain's vertical walls.
 In case of $Ra=10^5$ cells are moving apart from the center of the domain and it generates many such co-rotating cells which make the flow pattern more complex in nature.\\

  The thermal boundary layer is clearly visible on the left side of the domain starting from the hot wall of the domain depicted by fig.9, fig.11, fig.13.
  As Ra increased, beginning at $Ra=10^4$, it began to extend to the right side of the domain. When Ra value exceeded $Ra=10^4$, additional curves also began to appear close to the fins' walls. Two localised thermal zones, such as hot and cold zones, can be seen in the bottom left and top right of the domain for high values of Ra. When the value of $Ra$ was raised to $10^5$, they also started to move apart from one another toward the corner.
  The accompanying Matlab plotting, shown in fig.19, explains the impact of high value of Ra on heat transfer and change in thermal behaviour.
  
  The effect of Rayleigh number on Nu for a fixed value of $Pr(=1)$ can be observed from the table $5$. We can observe from table 5 that for $Pr=1$,
  $1 \leq Ra \leq 10^5$, Nu is an increasing function of $\phi,~0\% \leq \phi \leq 1\%$.
  Once more, we can see that Nu is an increasing function of Ra by placing $\phi$ within the range $0\% \leq phi \leq 1\%$ and $Pr=1$.
  
  The Matlab graphing provided by fig.14 and table $6$ allows one to observe the impact of the Prandtl number on Nu for a fixed value of $Ra(=10^4)$. 
  We can observe from table 6 that for $Ra=10^4$,
  $1\leq Pr \leq 10$, Nu is and increasing function of $\phi,~0\% \leq \phi \leq 1\%$.
  Again by fixing up $\phi$ within the range $0\% \leq \phi \leq 1\%$ and $Ra=10^4$ we can see that Nu is an increasing function of $Pr, ~1 \leq Pr \leq 10.$ The figures 15, 16, 17, and 18 also make it evident how Pr affects Nu value for a given value of Ra.
 For $Pr=1$, the streamlines in fig. 15 show the formation of two circulation eyes spaced equally apart from the vertical walls, whereas for $Pr=10$, fig. 17 shows the formation of just one main circulation eye located on the left side of the domain and the appearance of a multi-cellular pattern, which makes the flow more complex. Along with that presence of a very intensified multi-cellular  circulation zones are found in the left side of the domain near hot vertical wall. Figure 18 shows that an increase in Pr from $1$ to $10$ results in a highly noticeable hot zone in the bottom left corner of the domain and a cold thermal zone in the top right corner of the cavity. 
  \section{L shape geometry:}
 Here we have considered the following H shaped domain with $100X100$ grid size shown in Fig. 20. The numerical experiment is carried out with the following boundary conditions:\\
  Extreme left vertical wall is hot and is maintained at $T=1.$\\
  Extreme right vertical wall is cold and is maintained $T=0.$\\
  Rest of the walls of the domain are adiabatic.\\
  $U=V=0 $ throughout the boundary.\vspace{3mm}\\
 \begin{figure}[h]
  \centering
  \includegraphics[width=\textwidth]{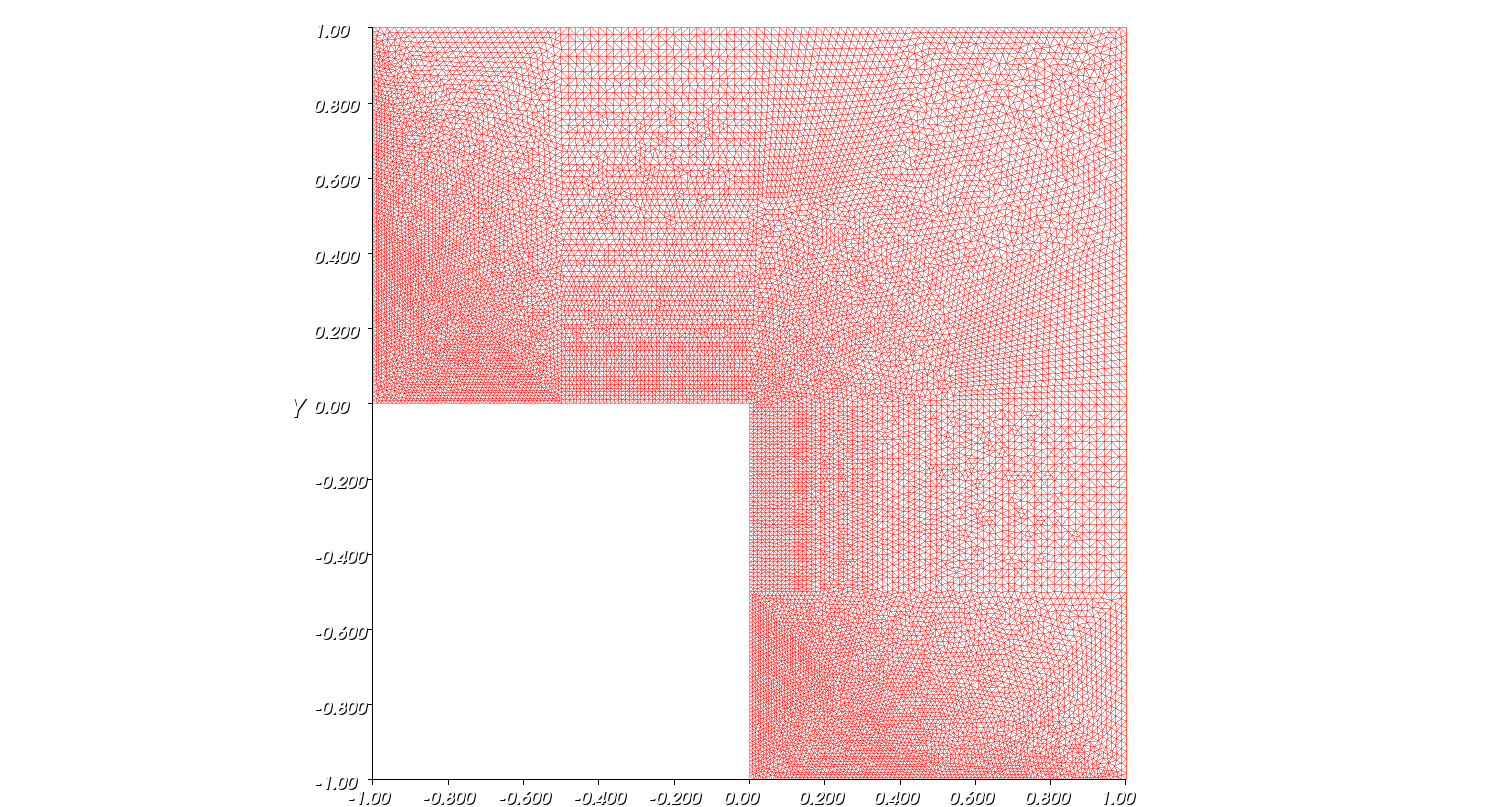}
    \caption{ L Shape Domain}
    \label{fig:height .25 unit,breadth 1 unit}
 \end{figure}
\subsection{ Comparison between the cases $Ra=10^5$ AND $Ra=10^6$ of HNF1,HNF2,HNF3:}
\begin{figure}[h!]
  \centering
  \begin{minipage}[b]{0.32\textwidth}
    \includegraphics[width=\textwidth]{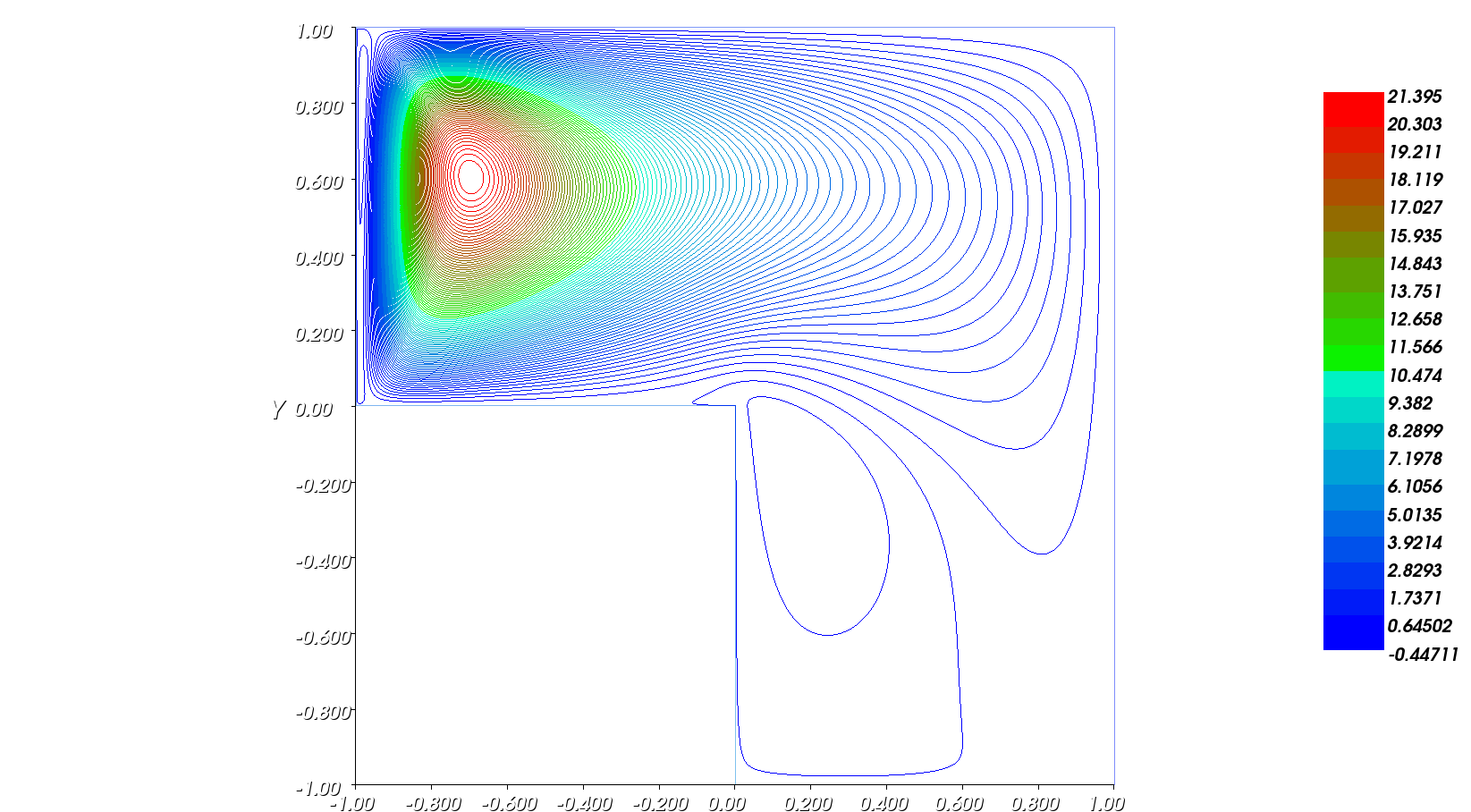}
    \caption{streamlines  }
  \end{minipage}
  \hfill
  \begin{minipage}[b]{0.32\textwidth}
    \includegraphics[width=\textwidth]{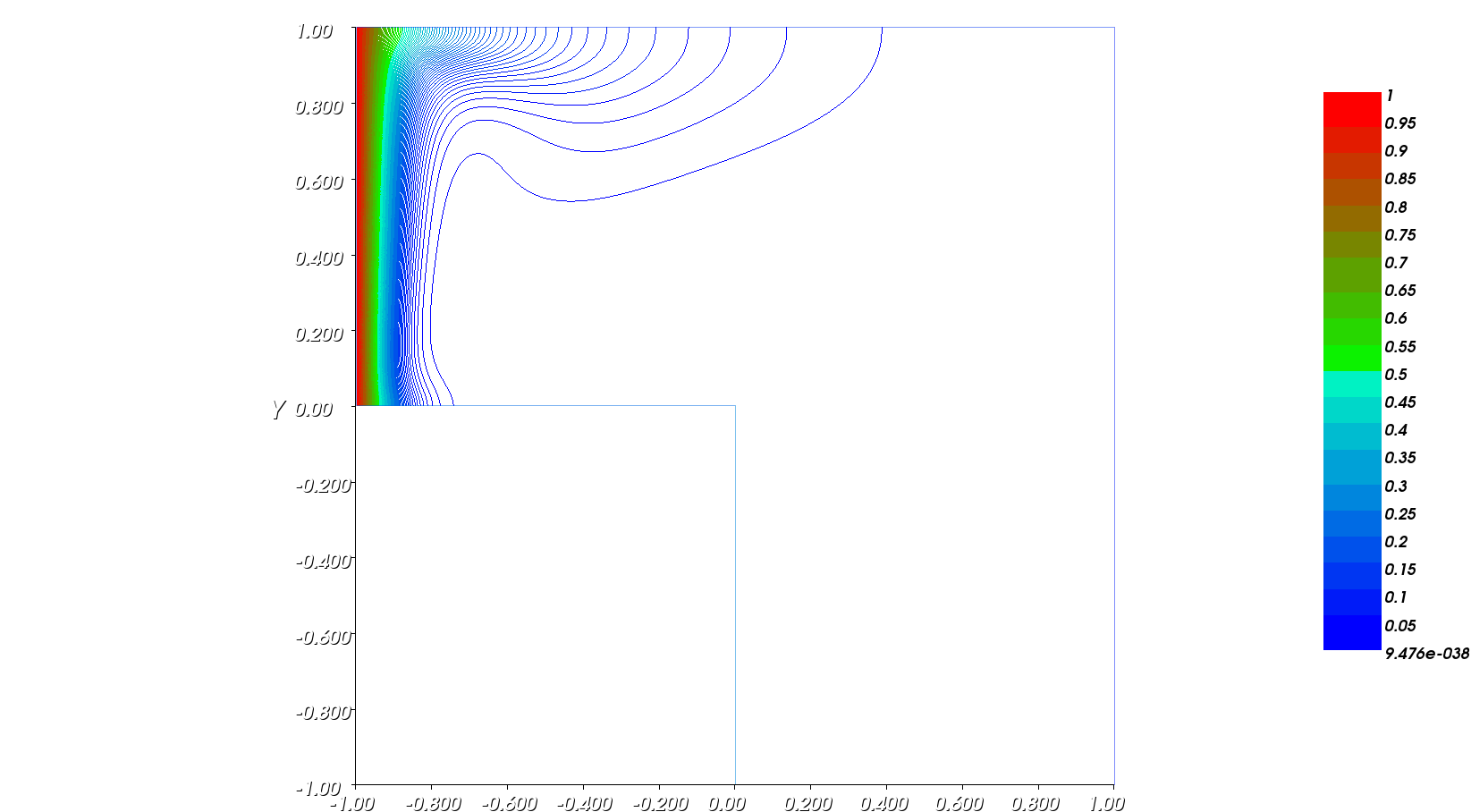}
    \caption{Isotherms }
  \end{minipage}
  \hfill
  \centering
  \begin{minipage}[b]{0.32\textwidth}
    \includegraphics[width=\textwidth]{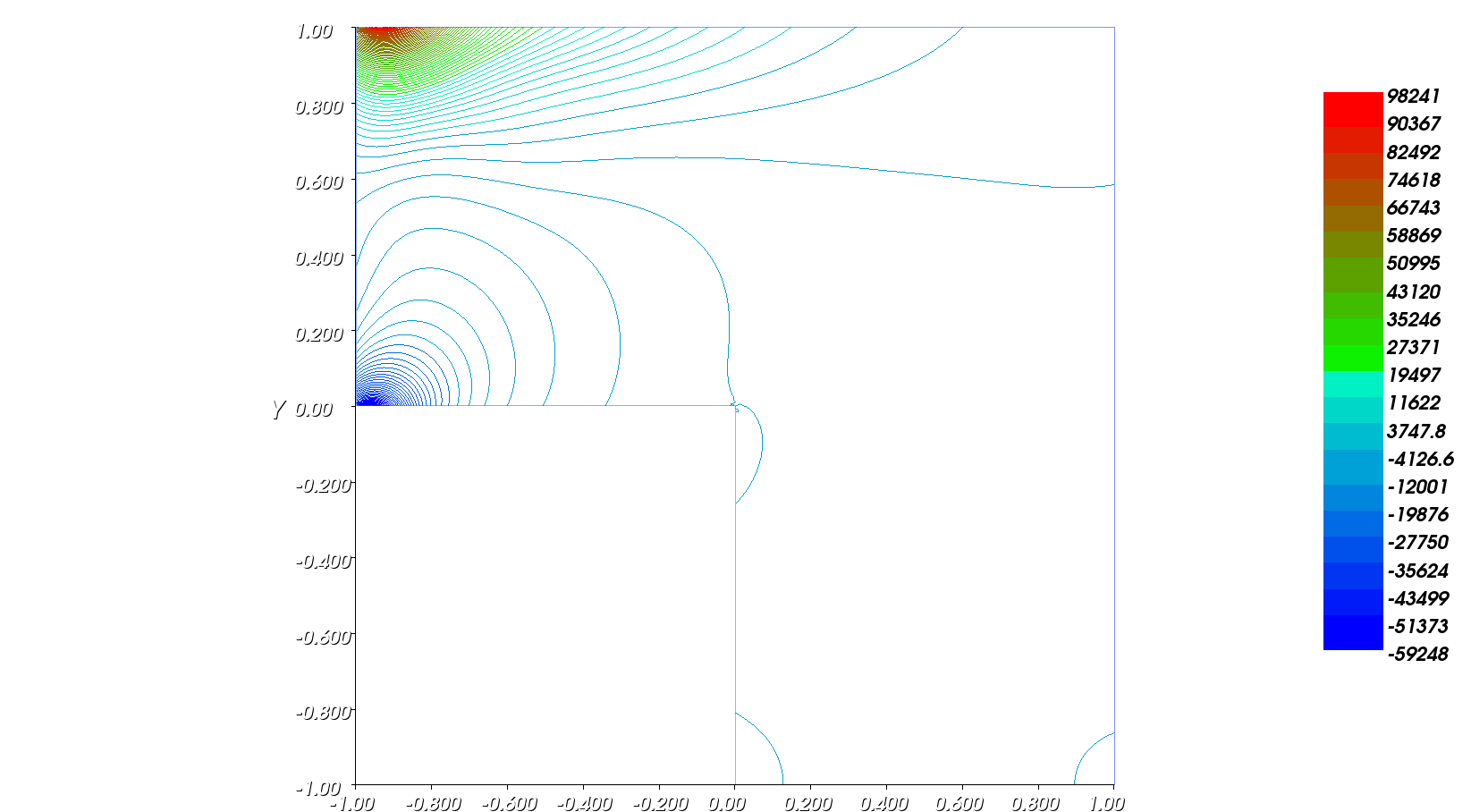}
    \caption{Pressure  }
  \end{minipage}
   ~~~~~~~~~~~~~~~~~~~~~~~~~~~~~~~~~~~~~~~~~~~~~~~~~~~~~~~~~~~~~~~~~~~~~~~~~~~~~~~~~~~~~~~Contours for $Ra=10^5$, HNF1
  \end{figure}

\begin{figure}[!h]
  \centering
  \begin{minipage}[b]{0.32\textwidth}
    \includegraphics[width=\textwidth]{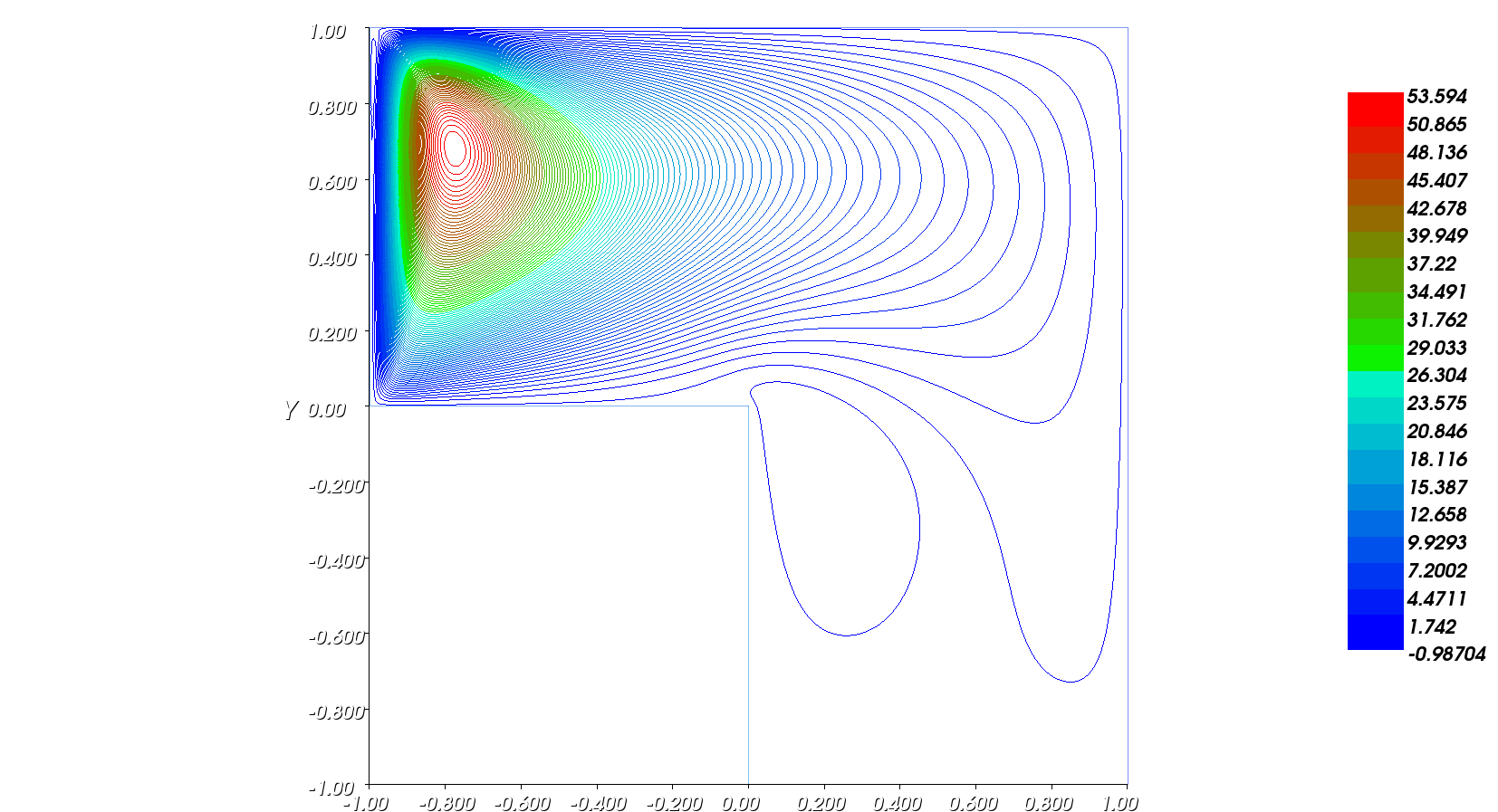}
    \caption{streamlines  }
  \end{minipage}
  \hfill
  \begin{minipage}[b]{0.32\textwidth}
    \includegraphics[width=\textwidth]{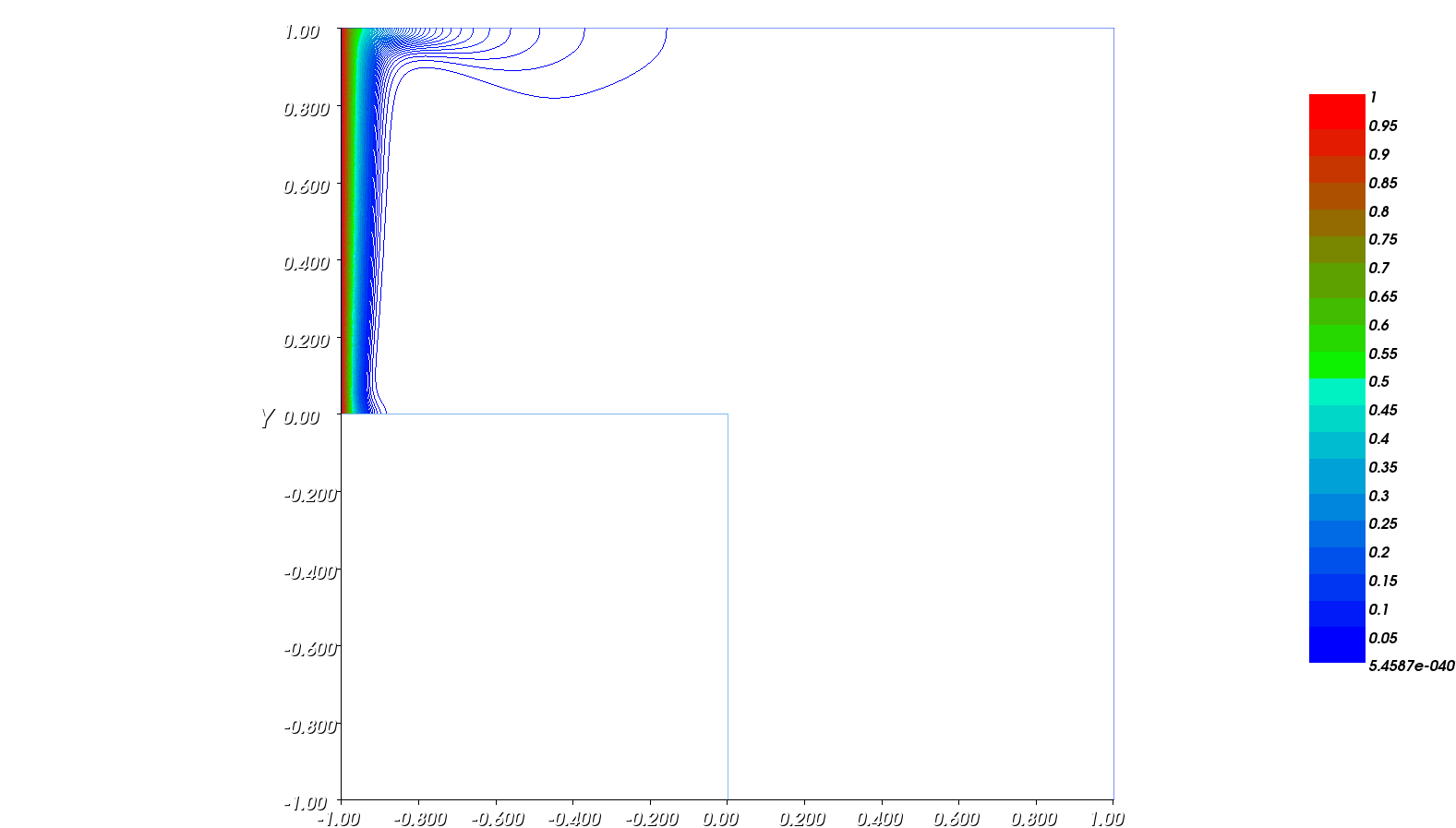}
    \caption{Isotherms }
  \end{minipage}
  \hfill
  \centering
  \begin{minipage}[b]{0.32\textwidth}
    \includegraphics[width=\textwidth]{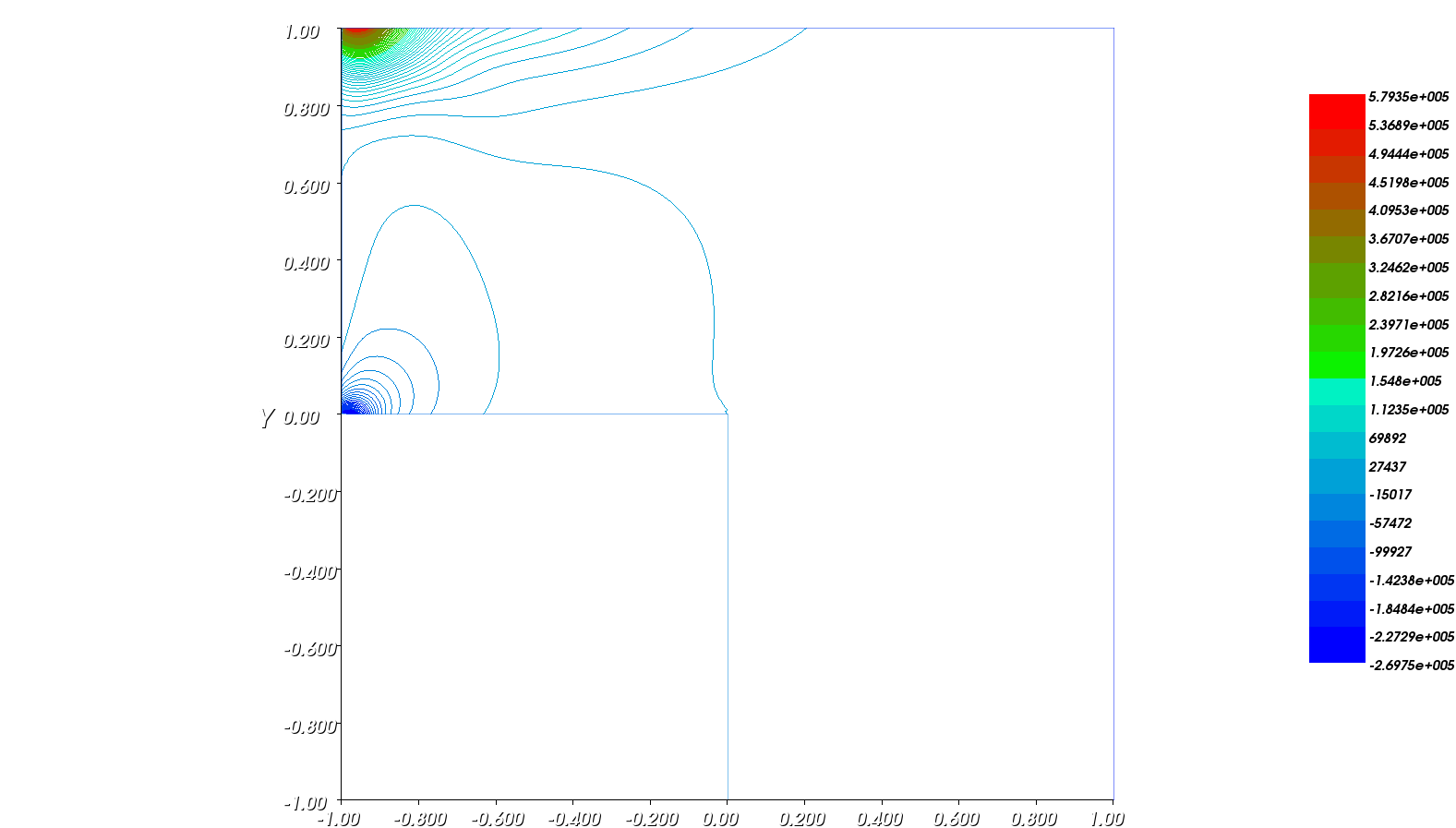}
    \caption{Pressure  }
  \end{minipage}
   ~~~~~~~~~~~~~~~~~~~~~~~~~~~~~~~~~~~~~~~~~~~~~~~~~~~~~~~~~~~~~~~~~~~~~~~~~~~~~~~~~~~~~~~Contours for $Ra=10^6$, HNF1.

  \end{figure}


\begin{figure}[!h]
  \centering
  \begin{minipage}[b]{0.32\textwidth}
    \includegraphics[width=\textwidth]{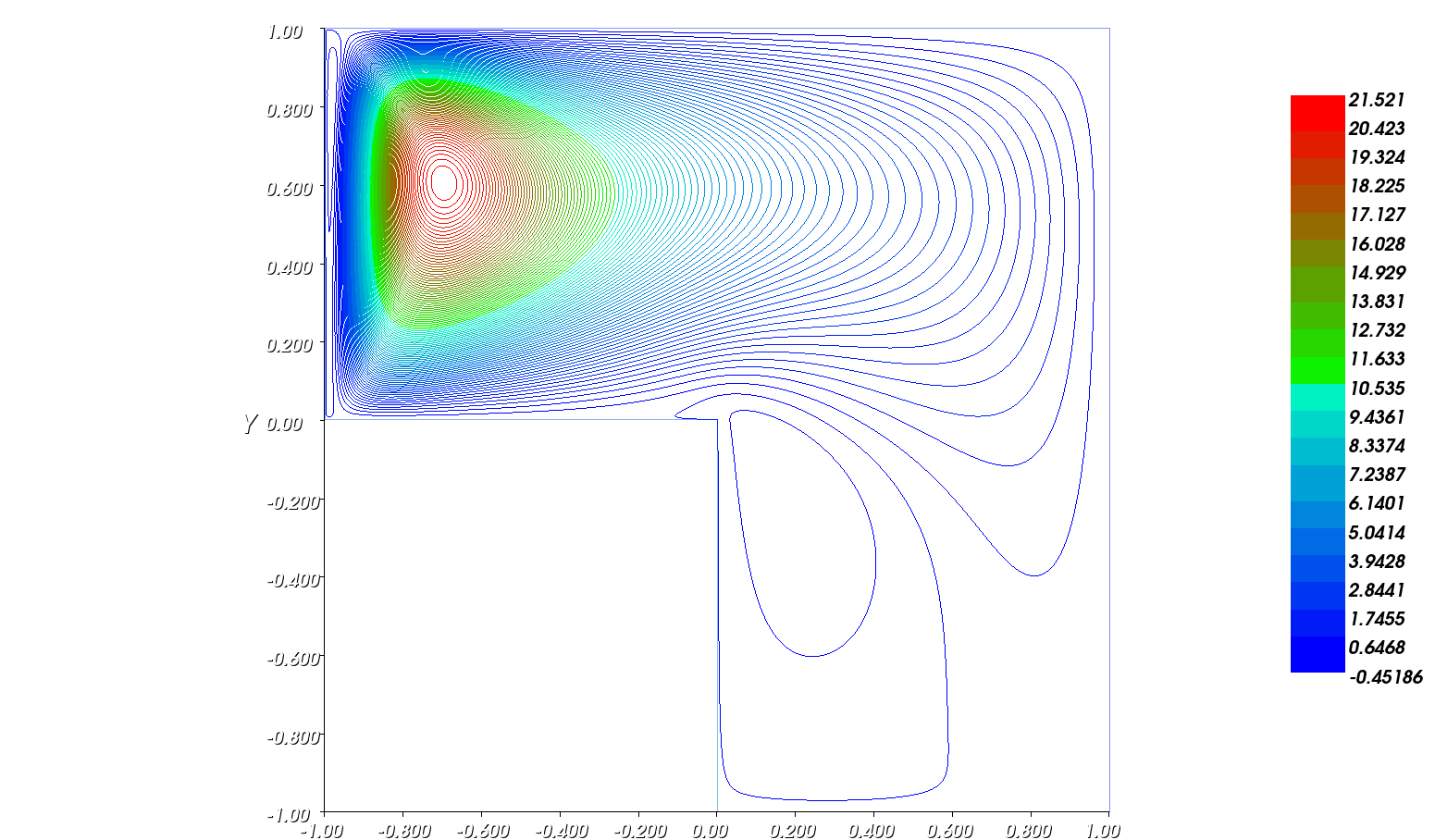}
    \caption{streamlines  }
  \end{minipage}
  \hfill
  \begin{minipage}[b]{0.32\textwidth}
    \includegraphics[width=\textwidth]{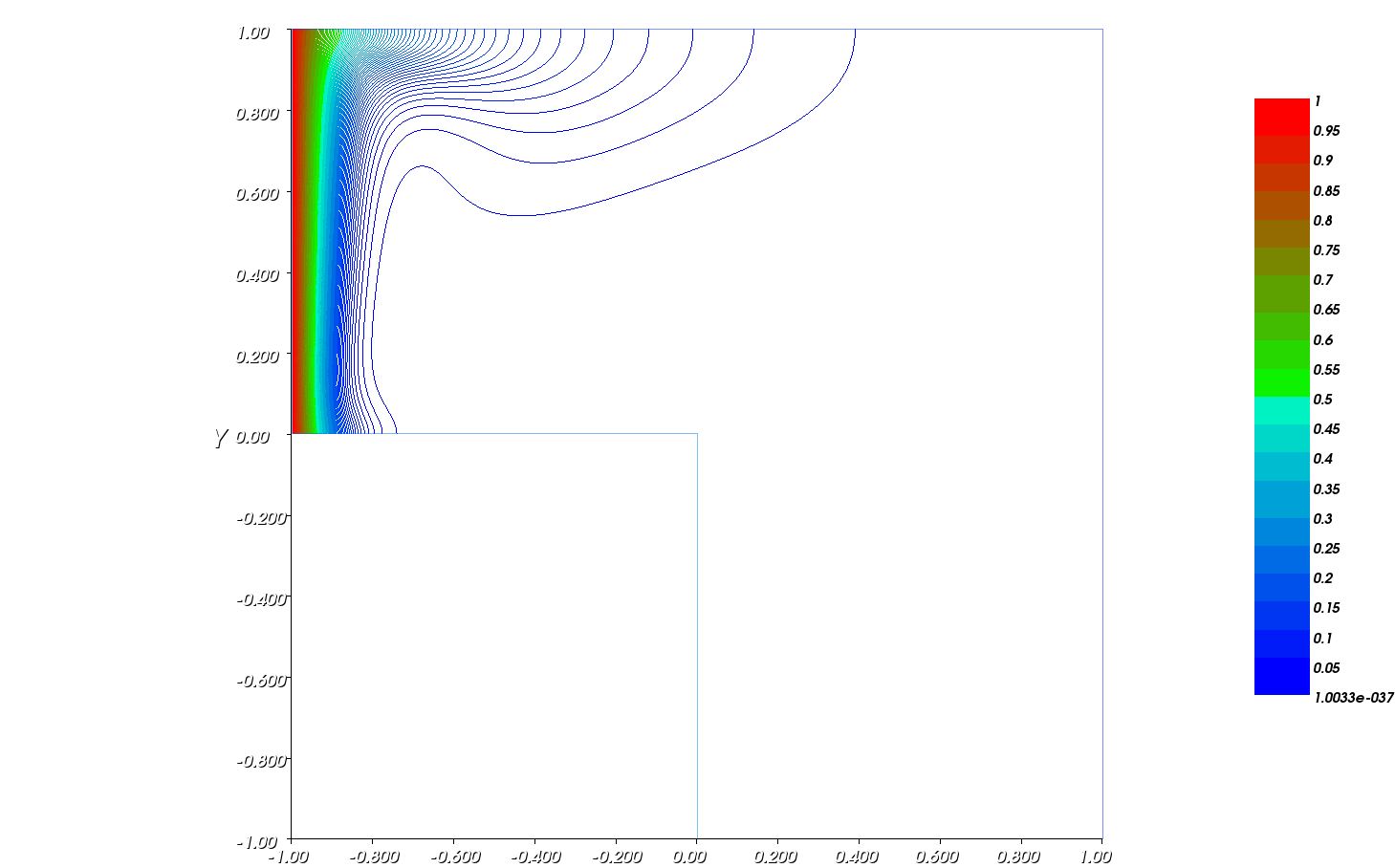}
    \caption{Isotherms }
  \end{minipage}
  \hfill
  \centering
  \begin{minipage}[b]{0.32\textwidth}
    \includegraphics[width=\textwidth]{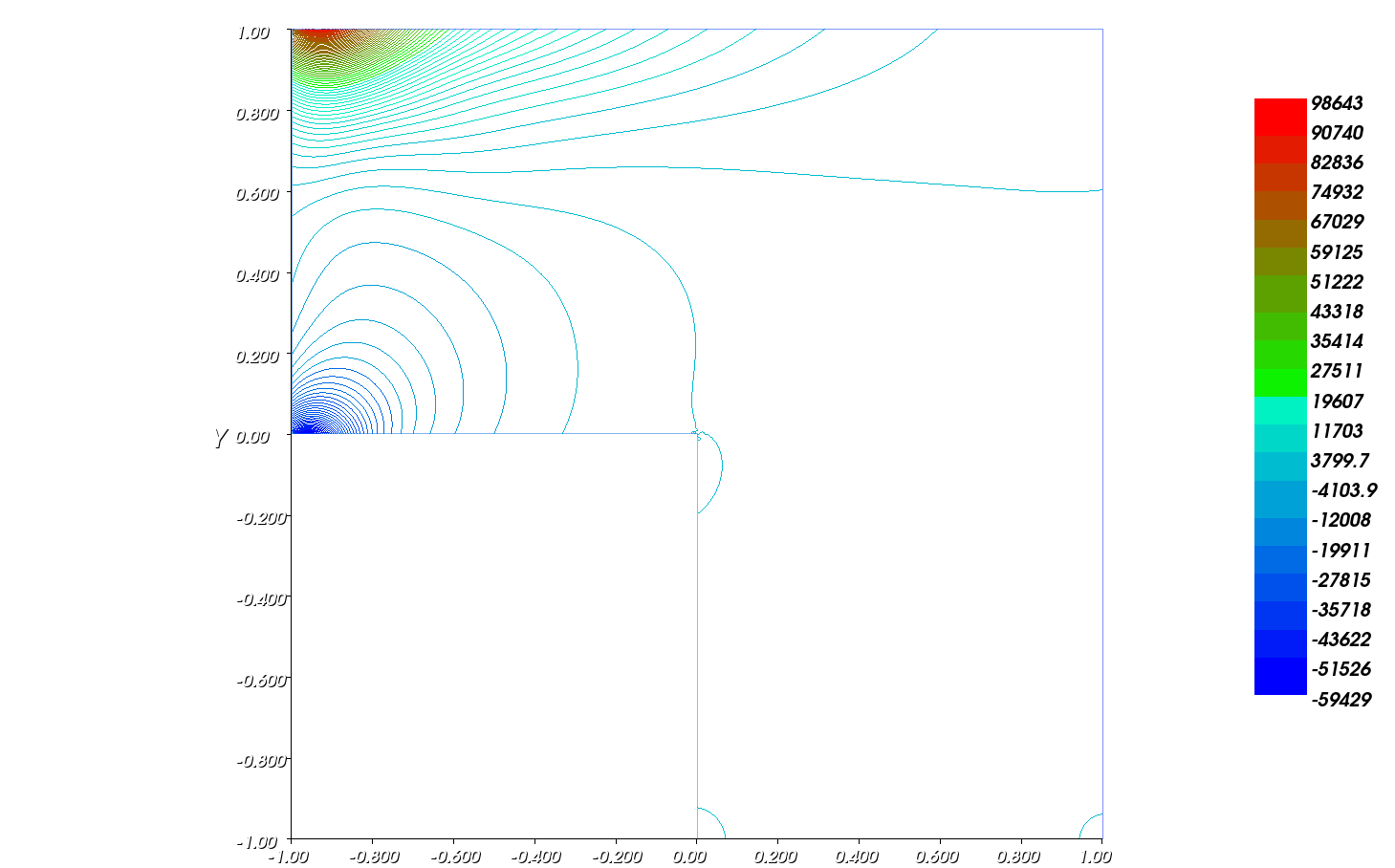}
    \caption{Pressure} 
  \end{minipage}
~~~~~~~~~~~~~~~~~~~~~~~~~~~~~~~~~~~~~~~~~~Contours for $Ra=10^5$, HNF2
  \end{figure}

\begin{figure}[!h]
  \centering
  \begin{minipage}[b]{0.32\textwidth}
    \includegraphics[width=\textwidth]{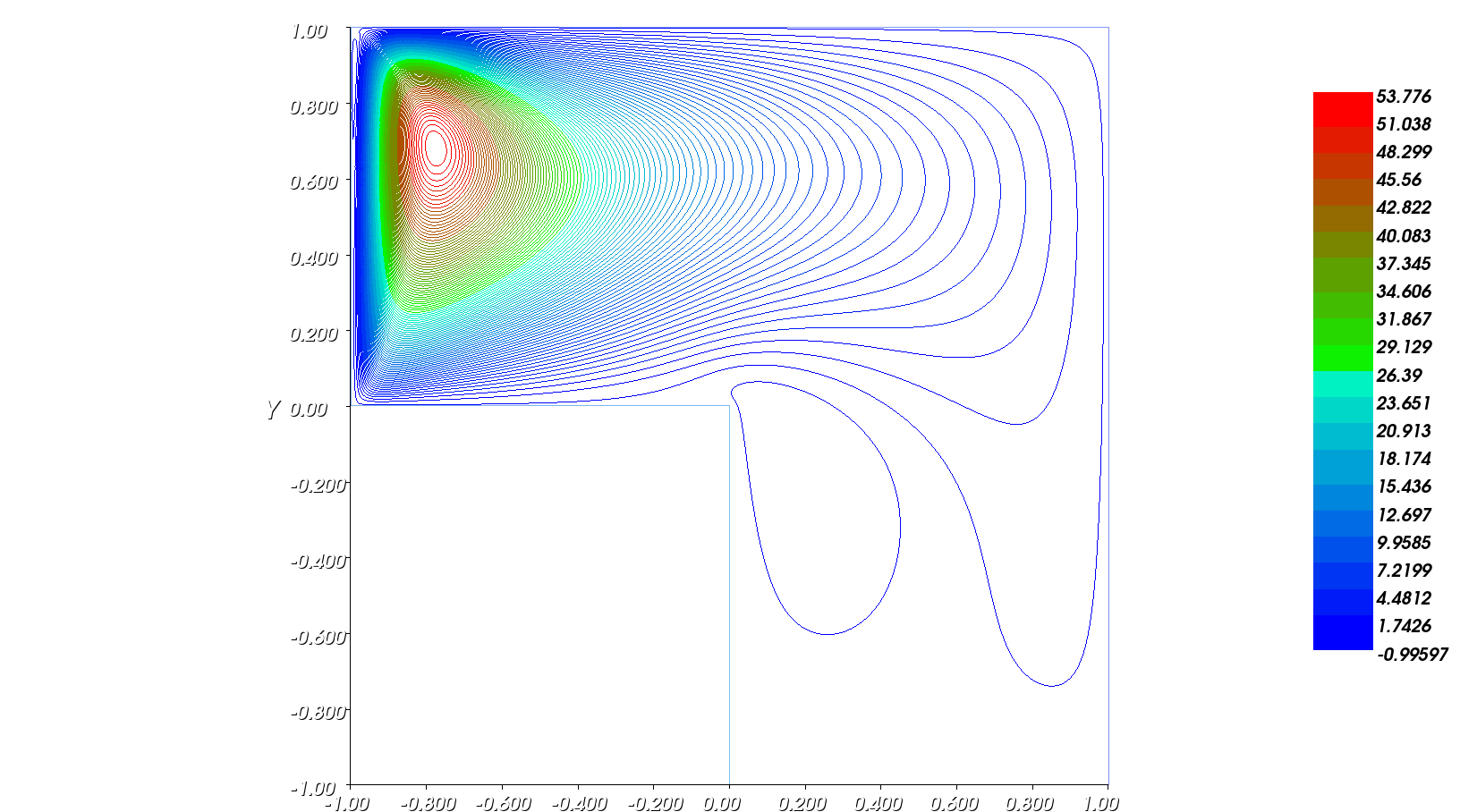}
    \caption{streamlines  }
  \end{minipage}
  \hfill
  \begin{minipage}[b]{0.32\textwidth}
    \includegraphics[width=\textwidth]{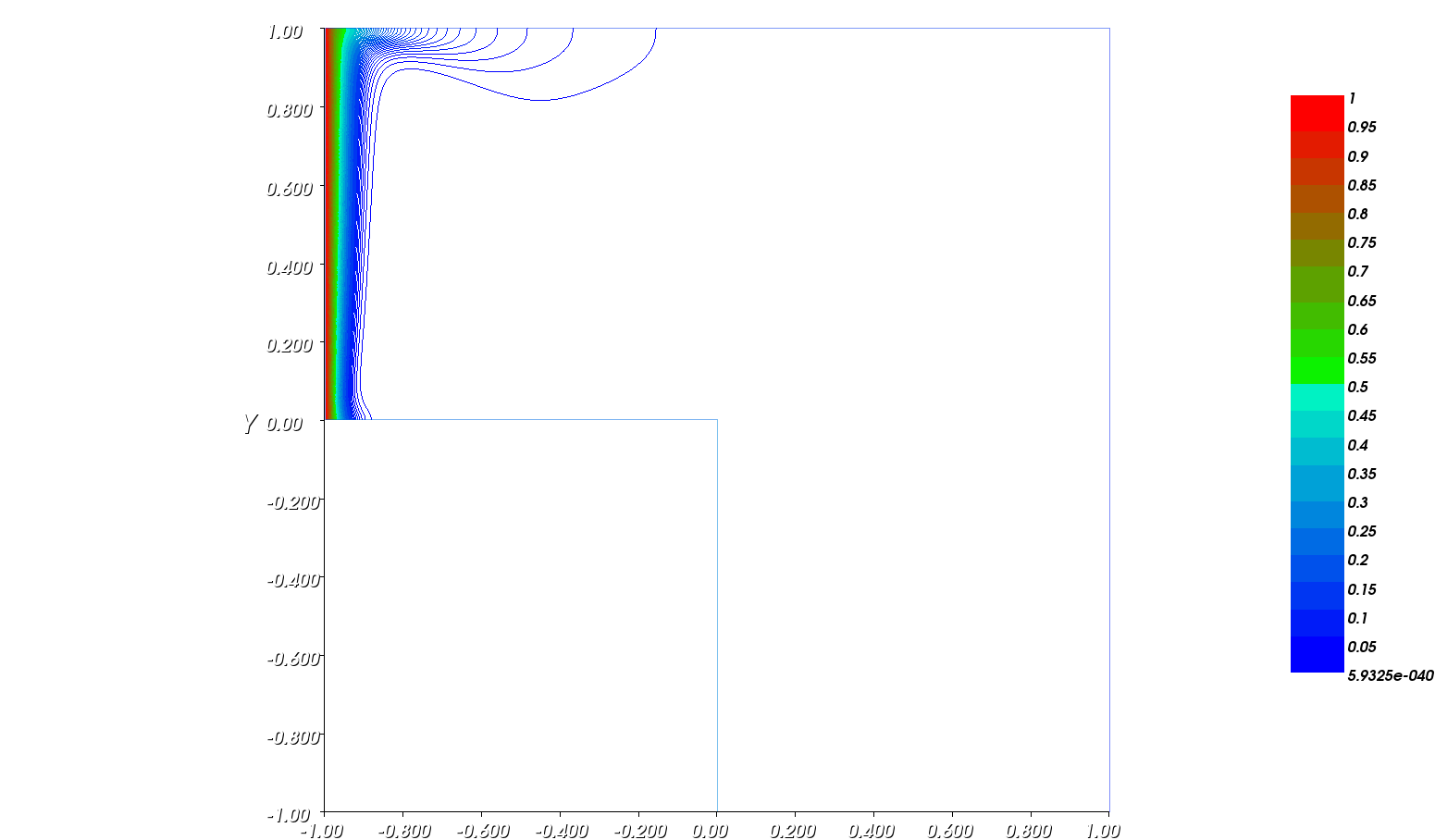}
    \caption{Isotherms }
  \end{minipage}
  \hfill
  \centering
  \begin{minipage}[b]{0.32\textwidth}
    \includegraphics[width=\textwidth]{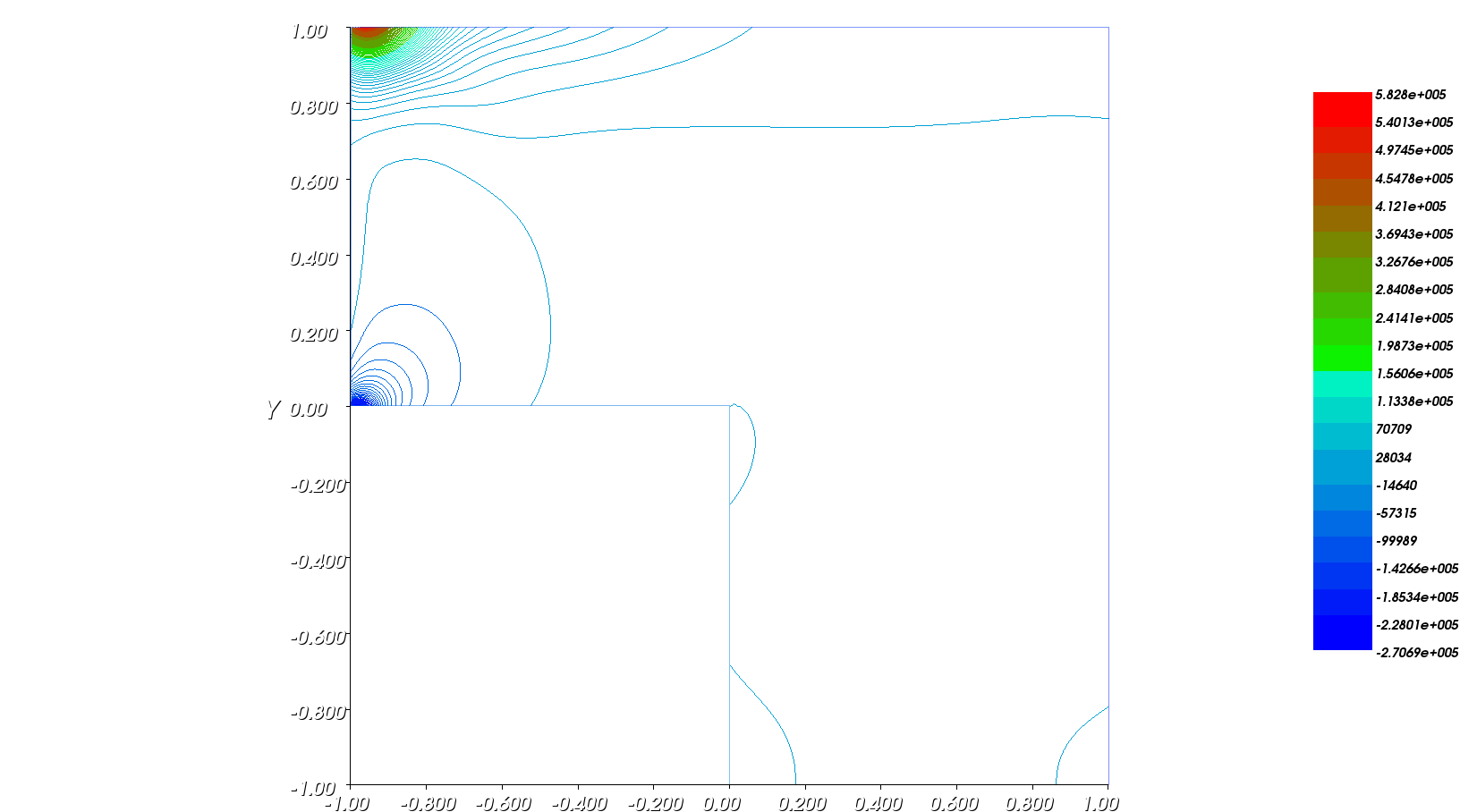}
    \caption{Pressure  }
  \end{minipage}
   ~~~~~~~~~~~~~~~~~~~~~~~~~~~~~~~~~~~~~~~~~~~~~~~~~~~~~~~~~~~~~~~~~~~~~~~~~~~~~~~~~~~~~~~Contours for $Ra=10^6$,HNF2

  \end{figure}
  

\begin{figure}[!h]
  \centering
  \begin{minipage}[b]{0.32\textwidth}
    \includegraphics[width=\textwidth]{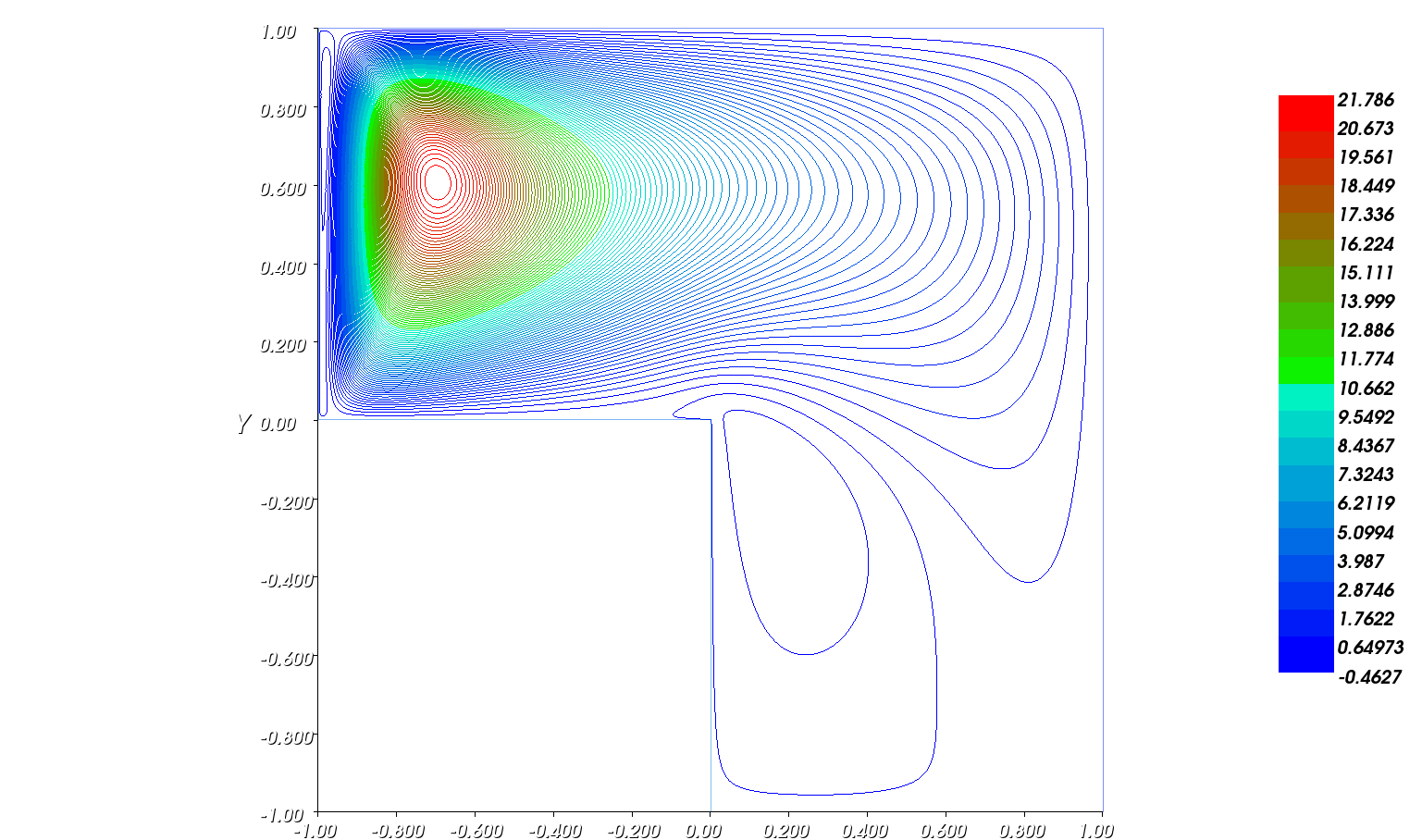}
    \caption{streamlines  }
  \end{minipage}
  \hfill
  \begin{minipage}[b]{0.32\textwidth}
    \includegraphics[width=\textwidth]{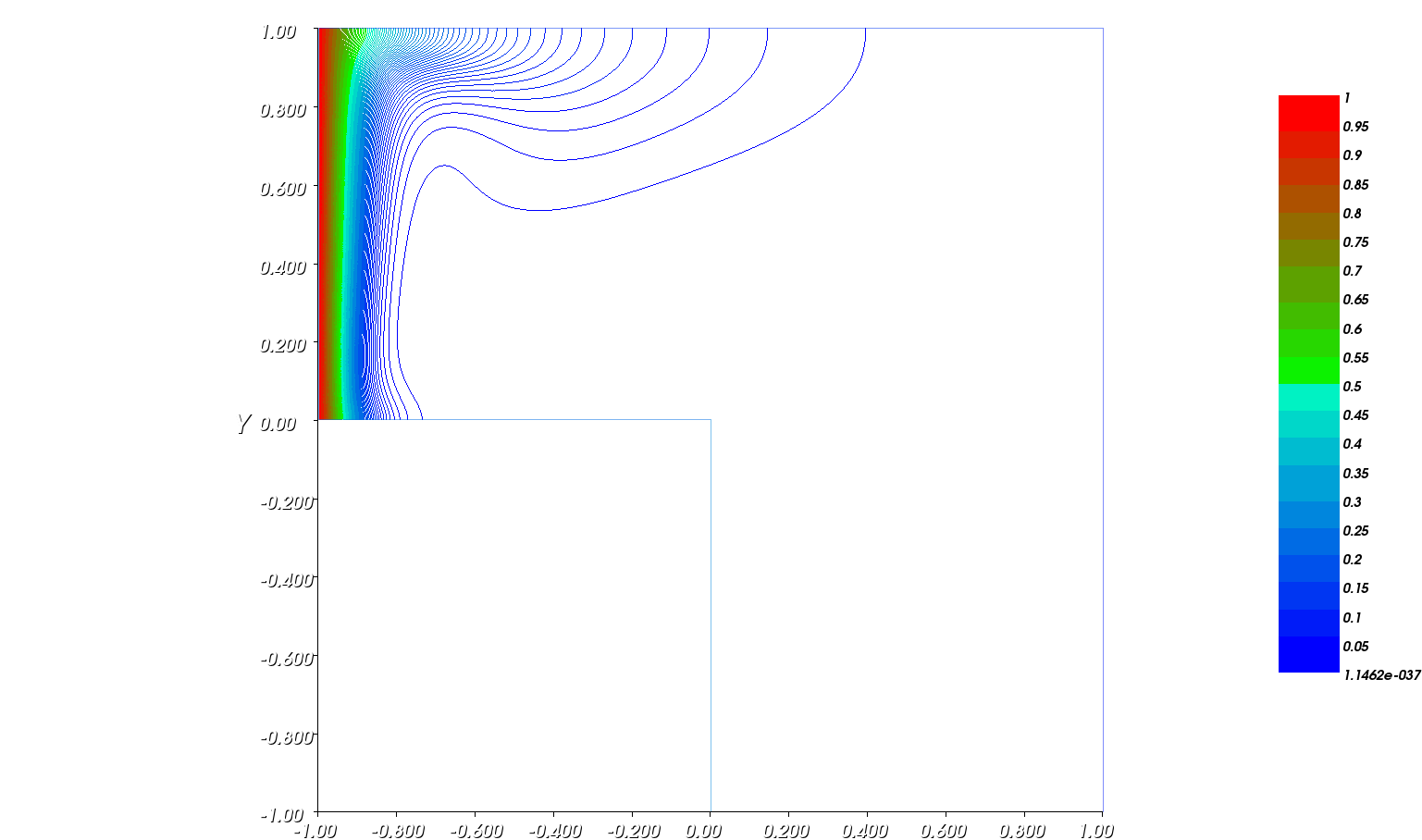}
    \caption{Isotherms }
  \end{minipage}
  \hfill
  \centering
  \begin{minipage}[b]{0.32\textwidth}
    \includegraphics[width=\textwidth]{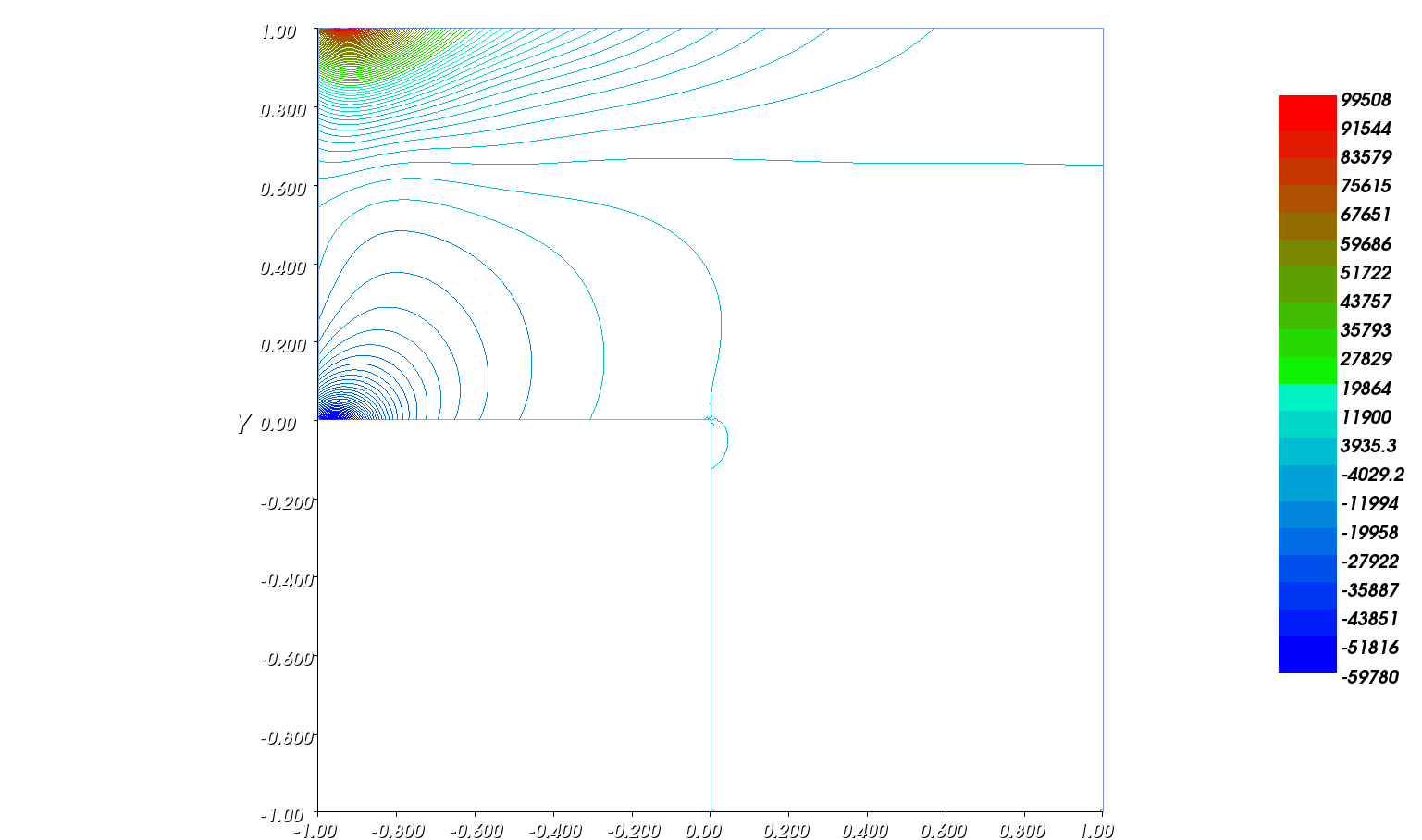}
    \caption{Pressure  }
  \end{minipage}
   ~~~~~~~~~~~~~~~~~~~~~~~~~~~~~~~~~~~~~~~~~~~~~~~~~~~~~~~~~~~~~~~~~~~~~~~~~~~~~~~~~~~~~~~Contours for $Ra=10^5$,HNF3

  \end{figure}

\begin{figure}[!h]
  \centering
  \begin{minipage}[b]{0.32\textwidth}
    \includegraphics[width=\textwidth]{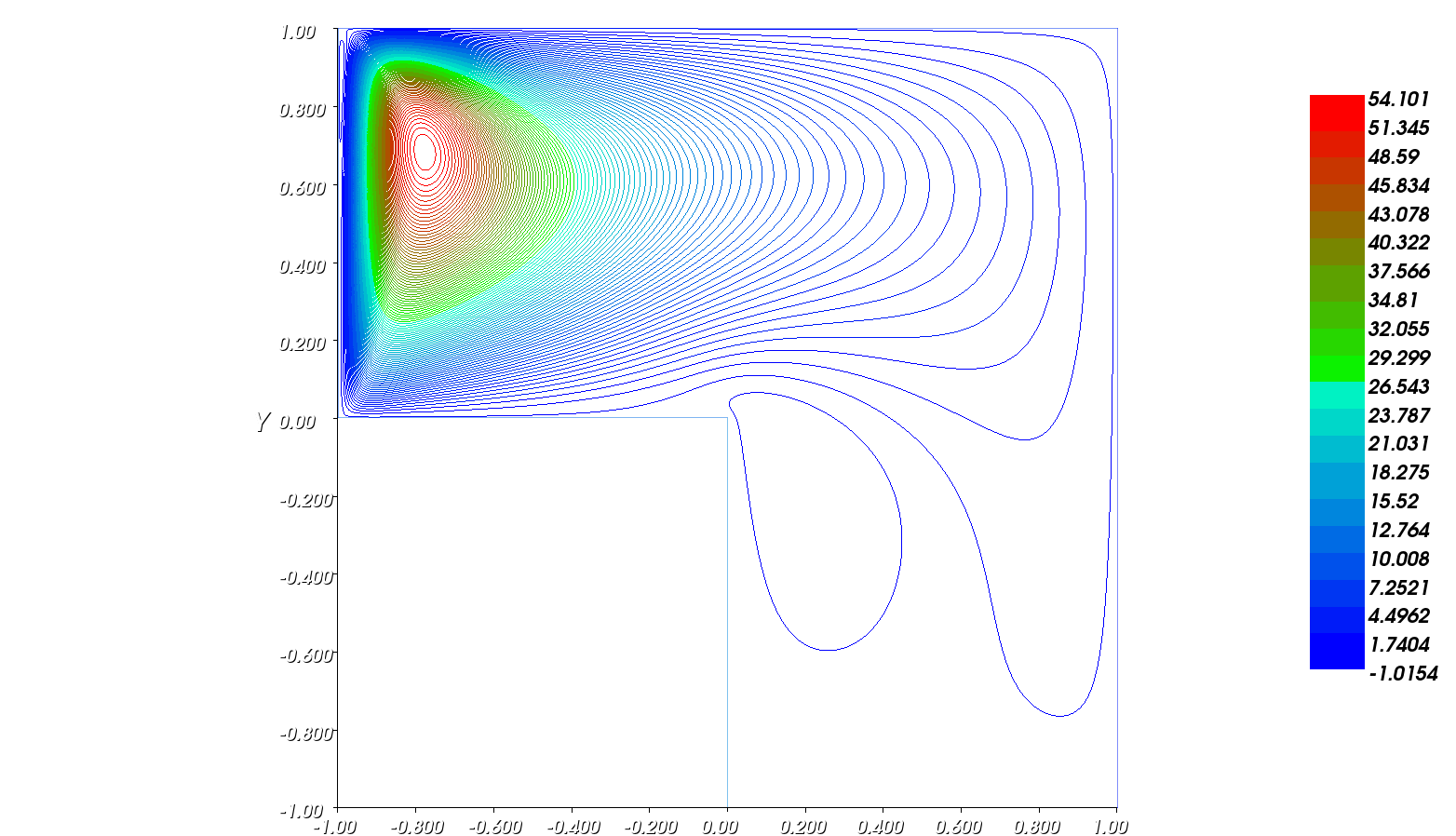}
    \caption{streamlines  }
  \end{minipage}
  \hfill
  \begin{minipage}[b]{0.32\textwidth}
    \includegraphics[width=\textwidth]{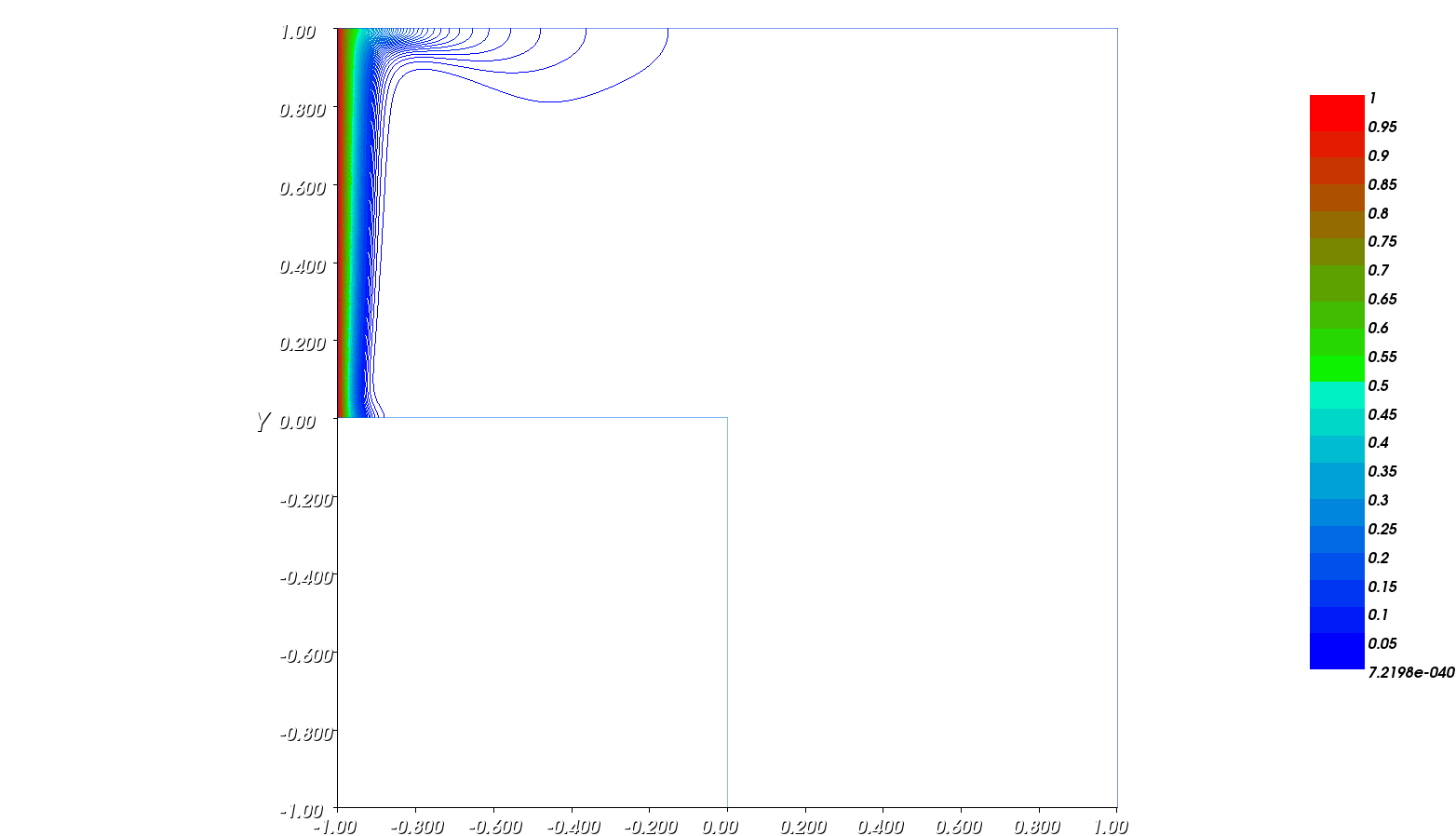}
    \caption{Isotherms }
  \end{minipage}
  \hfill
  \centering
  \begin{minipage}[b]{0.32\textwidth}
    \includegraphics[width=\textwidth]{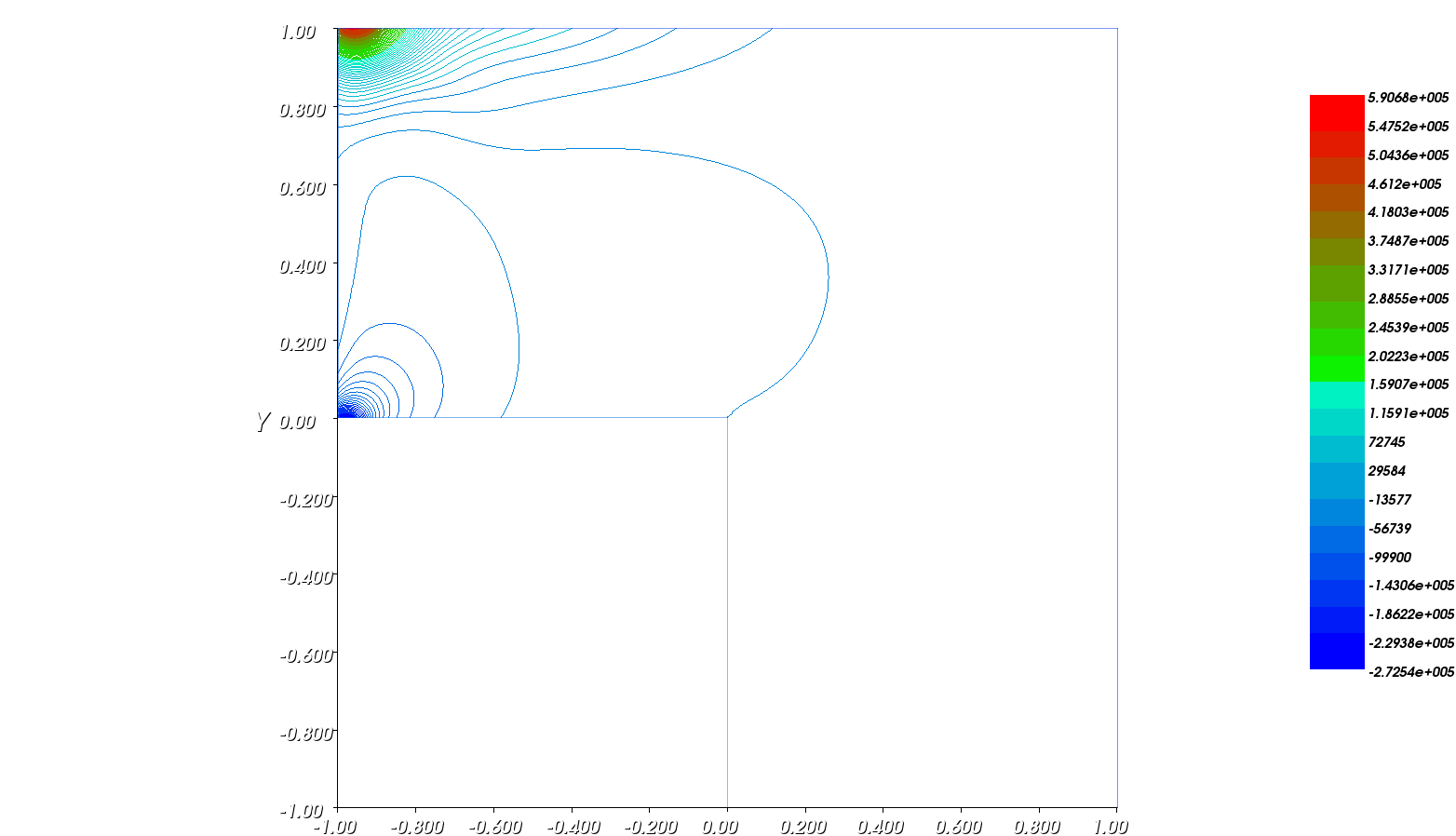}
    \caption{Pressure  }
  \end{minipage}
   ~~~~~~~~~~~~~~~~~~~~~~~~~~~~~~~~~~~~~~~~~~~~~~~~~~~~~~~~~~~~~~~~~~~~~~~~~~~~~~~~~~~~~~~Contours for $Ra=10^6$,HNF3

  \end{figure}
  \clearpage
  \newpage
  \subsection{Nusselt no. table:}
 Here we have considered three types of volume fractions of nano-particles as 0.1$\%$, 0.33$\%$, and 1.0$\%$. We have renamed the problems as HNF1, HNF2, HNF3 respectively.

\begin{table}[h]
	\centering
	\scalebox{5}{}
	\begin{tabular}{l|c|c|c|r}
		\hline
		\hline
		
			$Pr=10$,Ra  & HNF1 Nu & HNF2 Nu & HNF3 Nu \\
			\hline
			$Ra=10^5$  & 8.33827488 & 8.43228086 & 8.65425299\\
			\hline 
			$ Ra=3 \times 10^5$  & 11.8195846 &11.9610339 & 12.2837988 \\
			\hline 
			$Ra=5 \times 10^5$  &  13.8238296 & 13.9901545 & 14.3696117 \\
			\hline 
			$Ra=7 \times 10^5$ & 15.2999125& 15.4842286 & 15.9049525 \\
			\hline
			$Ra=10^6$ & 17.0135319 & 17.2187546 & 17.6868126\\
				\hline
		
	\end{tabular}
	\caption{Comparison of Nusselt no. for different volume fraction on L shape domain}
\end{table}

 \begin{figure}[h]
 \centering
  \includegraphics[width=90mm]{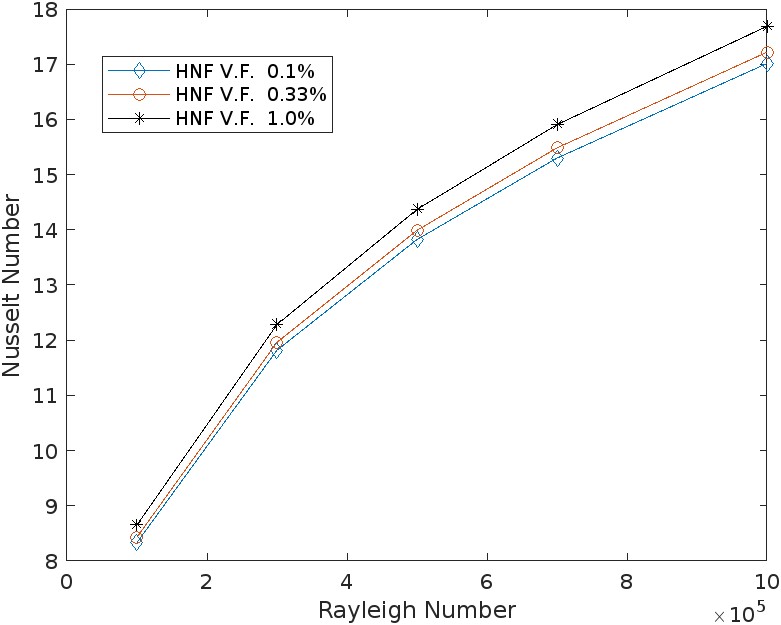} 
    \caption{Effect of Ra and volume fraction on Nu for fixed value of Pr}
    \label{fig:height .25 unit,breadth 1 unit}
  \end{figure} 
 
\subsection{Effect of source length on global Nusselt no.:}
Another important input parameter to influence heat transfer is source length and its effect is presented by the following table for  HNF1, HNF2 and HNF3 when $Pr=10$ and $Ra=10^5$.
 
\begin{table}[!h]
	\centering
	\scalebox{1.5}{}
	
	\begin{tabular}{|l|c|c|c|}
		\hline
		&\multicolumn{3}{c}{|Nusselt No.|} \\
		\cline{2-4}
		\newline
		Fluid & 1 unit  & 0.6 units & 0.2 units\\
		\hline
		HNF1 & 8.33488827   & 7.11625273 & 3.37201574 \\
		\hline
		HNF2 & 8.43228086  & 7.20183898 & 3.41463114  \\
		\hline
		HNF3   & 8.65328838 & 7.3967311 & 3.51205609\\
			\hline
	\end{tabular}
	\caption{Effect of Source length for $Pr = 10,Ra = 10^5$}
\end{table}

\subsection{Results and Discussions on L shape domain:}
Main center of circulation lies near the hot wall and formation of some more secondary vortices are noticed. 
Also formation of co-rotating circulation zones are discernible around the main eye of the circulation.\\
With the isotherm contours one can notice that when Ra rises from $10^5$ to $10^6$ then the number of curves decrease and they get concentrated along the hot wall. A clear thermal boundary layer is visible. Isotherm curves start from the left hot boundary which reach up to the top horizontal wall.\\
The high pressure region is situated in the top left corner of the domain whereas the low pressure region is noted just below that region close to the origin.\\ An increase in $Ra$ value from $10^5$ to $10^6$ causes the centralization of the low and high pressure region towards the middle left and top left corners respectively. Additionally, appearance of more number of curves is noticed when thermal buoyancy force is augmented from  $Ra=10^5$ to $10^6.$
On observing the data in the table 7 one can see that for $0.1\% \leq \phi \leq 1\%$, Nu is an increasing function of $\phi$ for each fixed value of Ra. Also the increase in  thermal buoyancy forces resulted the enhancement of heat transfer, for a fixed value of Pr, as is once more evident from the Matlab plotting provided by fig. 39. Effect of heat source on heat transfer is clear from the table 8. For a fixed value of $\phi$, $0.1\% \leq \phi \leq 1\%$, we can observe a fall in the values of Nu with the decreasing of source length whereas for each source length, Nu is an increasing function of $\phi.$
\clearpage
\newpage
\section{Conclusion:}
By using finite element method we have analysed convective flow and heat transfer in  H and L shaped complex domains. A detailed convergence analysis through apriori error estimation is established for the finite element method. We have carried out the detailed of grid and code validation tests before proceeding for the numerical simulation.  The major inferences from the study are summarized below:
\begin{itemize}
    
\item{ Apriori error estimation is derived in $H^{1}$ norm for the proposed  Galerkin Finite Element  Approach.}\\

\item{ In case of both the H and L shaped domains a boost in thermal Buoyancy forces measured in terms of Ra results in increase in Nu.}\\

\item{ Whether it is H shaped or L shaped domain,  for a fixed value of $Ra(=10^4)$ an augmentation in Pr value causes an enhancement in convective fluxes i.e. Nu value.}\\

\item{ Our numerical experiment discloses that Nu is an increasing function of volume fraction $\phi.$}\\

\item{ In the case of partially heated side walls of L shaped cavity, a fall in global heat flux (Nu) values is observed with the decreasing of source length. }

\end{itemize}



\end{document}